\documentclass[preprints,article,accept,pdftex,moreauthors]{mdpi} 

\newcommand{\ncm}{\newcommand}
\ncm{\R}{{\mathbb{R}}}
\ncm\re[1]{(\ref{#1})}
\ncm\lb[1]{\label{#1}}
\ncm{\hI}{\hat{I}}
\ncm{\hQ}{\hat{Q}}
\ncm{\cP}{{\cal P}}
\ncm{\la}{\lambda}
\ncm{\Lto}{\, \stackrel{\cal L} \longrightarrow}
\ncm{\bfm}[1]{\mbox{\boldmath $#1$}}
\ncm{\bm}{\bfm{m}}
\ncm{\Si}{\Sigma}
\ncm{\bSi}{\bfm{\Si}}
\ncm{\si}{\sigma}
\ncm{\hxi}{\hat{\xi}}
\ncm{\hth}{\hat{\theta}}
\ncm{\cS}{{\cal S}}
\ncm{\cA}{{\cal A}}
\ncm{\cX}{{\cal X}}
\ncm{\ga}{\gamma}
\ncm{\bPi}{\bfm{\Pi}}
\ncm{\pto}{\, \stackrel{p} \longrightarrow}

\newtheorem{prop}{Proposition}

{\theoremstyle{definition}
\newtheorem{exa}{\bf Example}
\newtheorem{rem}{\bf Remark}
}

\firstpage{1} 
\pubvolume{1}
\issuenum{1}
\articlenumber{0}
\pubyear{2022}
\copyrightyear{2022}
\datereceived{} 
\dateaccepted{} 
\datepublished{} 
\hreflink{https://doi.org/} 
\pdfoutput=1

\Title{Assessing, testing and estimating the amount of fine-tuning by means of active information}

\TitleCitation{Fine-tuning and active information}

\Author{Daniel Andr\'es D\'iaz-Pach\'on $^{1}$\orcidA{} and Ola H\"ossjer $^{2,}$*\orcidB{}}

\AuthorNames{Daniel Andr\'es D\'iaz-Pach\'on and Ola H\"ossjer}

\AuthorCitation{D\'iaz-Pach\'on, D. A.; H\"ossjer, O.}

\address{%
$^{1}$ \quad Division of Biostatistics, University of Miami; DDiaz3@miami.edu\\
$^{2}$ \quad Department of Mathematics, Stockholm University; ola@math.su.se}

\corres{Correspondence: ola@math.su.se; Tel.:  +46-706721218 (O.H.)}

\abstract{A general framework is introduced to estimate how much external information has been infused into a search algorithm, the so-called active information. This is rephrased as a test of fine-tuning, where tuning corresponds to the amount of pre-specified knowledge that the algorithm makes use of in order to reach a certain target. A function $f$ quantifies specificity for each possible outcome $x$ of a search, so that the target of the algorithm is a set of highly specified states, whereas fine-tuning occurs if it is much more likely for the algorithm to reach the target than by chance. The distribution of a random outcome $X$ of the algorithm involves a parameter $\theta$ that quantifies how much background information that has been infused. A simple choice of this parameter is to use $\theta f$ in order to exponentially tilt the distribution of the outcome of the search algorithm under the null distribution of no tuning, so that an exponential family of distributions is obtained. Such algorithms are obtained by iterating a Metropolis-Hastings type of Markov chain, and this makes it possible to compute the their active information under equilibrium and non-equilibrium of the Markov chain, with or without stopping when the targeted set of fine-tuned states has been reached. Other choices of tuning parameters $\theta$ are discussed as well. Nonparametric and parametric estimators of active information and tests of fine-tuning are developed when repeated and independent outcomes of the algorithm are available. The theory is illustrated with examples from cosmology, student learning, reinforcement learning, a Moran type model of population genetics, and evolutionary programming.}

\keyword{Active information; exponential tilting; fine-tuning; functional information; large deviations; Markov chains; Metropolis-Hastings; Moran model; statistical estimation and testing} 

\begin{document}

\section{Introduction}


When G\"odel published his incompleteness theorems \citep{Godel1931}, there was a commotion in the mathematical world from which it has neither yet recovered nor fully assimilated the consequences \cite{Hofstadter1999}. Hilbert's program to base mathematics on a finite set of axioms had earlier been pursued by Alfred North Whitehead and Bertrand Russell \cite{WhiteheadRussell1927}. But this approach turned out to be wrong when G\"odel proved that no finite set of axioms in a formal system can prove all its true statements, including its own consistency. In similar but lesser scale, when David Wolpert and William MacReady published their No Free Lunch Theorems (NFLTs, \cite{WolpertMacReady1995,WolpertMacReady1997}), there was disquiet in the community because these results imply that there is no one-size-fit-all algorithm that can do well in all searches \cite{Wolpert2021}, so that a ``theory of everything'' is not possible in machine learning. Wolpert and MacReady concluded that it was necessary to incorporate ``problem-specific knowledge into the behavior of the algorithm'' \cite{WolpertMacReady1997}. Thus active information (actinfo) was introduced in order to measure the amount of information carried by such problem-specific knowledge \cite{DembskiMarks2009a,DembskiMarks2009b}. More specifically, the NFLTs say that no search works better on average than a blind search, i.e., a search according to a uniform distribution. Accordingly, actinfo is defined as
\begin{align}\label{I+}
	I^+ = \log \frac{P(A)}{P_0(A)},
\end{align}
where $A\subset\Omega$ is the non-empty target of the search algorithm, a subset of the finite sample space $\Omega$, and $P_0$ is a uniform probability measure ($P_0(A)=|A|/|\Omega|$). $P$ must be seen here as the probability measure induced by the problem-specific knowledge of the researcher, whereas $P_0$ is the underlying distri\-bu\-tion assumed in the NFLTs. It corresponds to  absence of problem specific knowledge, in accordance with Bernoulli's Principle of Insufficient Reason (PoIR). An equivalent characterization of actinfo is the reduction of functional information 
\begin{equation}
I^+ = I_{f_0} - I_f = - \log P_0(A) - (-\log P(A))
\label{I+2}
\end{equation}
between algorithms that do not and do make use background knowledge. The name functional information was introduced by Szostak and collaborators \cite{ HazenEtAl2007,Szostak2003}. It refers to applications where $A$ corresponds to all outcomes of an algorithm that are functional according to some criterion. Then $I_{f_0}$ and $I_f$ are the self-information (measured in nats) of the event that an algorithm $X$ produces a functional outcome, given that it was generated under $P_0$ and $P$, respectively.    

Suppose we do not know whether problem specific knowledge has been used or not when the random search $X\in\Omega$ was generated. This corresponds to a hypothesis testing problem 
\begin{equation}
\begin{array}{ll}
H_0: X \sim P_0,\\
H_1: X \sim P,
\end{array}
\label{H0H1}
\end{equation}
where data is generated from distributions $P_0$ and $P$ under the null and alternative hypotheses $H_0$ and $H_1$, respectively. It follows from (\ref{I+}) that $I^+$ is the log likelihood ratio when testing $H_0$ against $H_1$, if data is censored so that only $X\in A$ is known. 

When the sample space $\Omega$ is finite or a bounded subset of a Euclidean space, the PoIR can be motivated by the fact that the uniform distribution maximizes Shannon entropy, since thereby it maximizes ignorance about the outcome of $X$. However, the uniform distribution is not a feasible choice of $P_0$ for unbounded samples spaces. For this reason actinfo has been generalized to deal with unbounded spaces \cite{DiazMarks2020a}, by choosing $P_0$ to maximize Shannon entropy under side constraints $\xi$, such as existence of various moments. This gives rise to a family of null distributions $P_0=P_{0\xi}$, with a $\xi$ a nuisance parameter that has to be estimated or controlled for in order to estimate or give bounds for the active information. 

Actinfo has also been used for mode detection in unsupervised learning, among other applications \cite{DiazEtAl2019,LiuEtAl2022}. Based on previous work by Monta\~nez \cite{Montanez2017, Montanez2018}, actinfo has been used in the past for hypothesis testing \cite{DiazSaenzRao2020}. More specifically, \cite{DiazSaenzRao2020} regards $P$ as a random measure so that $I^+$ is random as well, and finds expressions for the tail probability of $I^+$. 

\subsection{Fine-tuning}

Fine-tuning (FT) was introduced by Carter in physics and cosmology \cite{Carter1974}. Accor\-ding to FT, the constants in the laws of nature and/or the boundary conditions in the standard models of physics must belong to intervals of low probability in order for life to exist. Since its inception, FT has generated a great deal of fascination, seen in multiple divulgation books (e.g., \cite{BarrowTipler1988, Davies1982, LewisBarnes2016,Rees2000}) and scientific articles (e.g., \cite{Adams2019, Barnes2012, TegmarkRees1998, TegmarkEtAl2006}). For a given constant of nature $X$, the connection between FT and active information can be described in three steps: 
\begin{enumerate}
	\item[(i)] Establishing the life-permitting interval (LPI) $A$ that allows the existence of life for the constant, with $\Omega=(0,\infty)=\R^+$ or $\Omega=\R$ the range of values this constant could possibly take, including those that do not permit life. 
	\item[(ii)] Determining the probability $P_0(A)$ of such a LPI. If $P_0=P_{0\xi}$ contains unknown parameters $\xi$, find an upper bound 
\begin{equation}
P_{0\text{max}}(A)=\max_\xi P_{0\xi}(A)
\label{P0max}
\end{equation}
of $P_0(A)$.
	\item[(iii)] Suppose $H_1$ corresponds to an agent who uses background knowledge of what is required for life to exist in order to bring about a constant of nature $X$ that with certainty permits life ($P(A)=1$). The active information $I^+ = I_{f_0} = -\log P_0(A)$ is then a measure of how much background knowledge this agent infused. Following \cite{DiazHossjerMarks2021, DiazHossjerMarks2022}, we conclude that $X$ is finely tuned when the lower bound $-\log P_{0\text{max}}(A)$ of $I^+=I_{f_0}$ is large enough. That is, FT corresponds to infusing a high degree of background knowledge into a problem.   
\end{enumerate}

Fine-tuning has also been used in biology. Dingjan and Futerman explored it for cell membranes \cite{DingjanFuterman2021a, DingjanFuterman2021b}, whereas Thorvaldsen and H\"ossjer \cite{ThorvaldsenHossjer2020} formalized it for a large class of biological models. According the \cite{ThorvaldsenHossjer2020}, a system is fine-tuned if it (a) has an independent specifica\-tion, and (b) is very unlikely to occur by chance.

\subsection{The present article}

In this article actinfo will not only be used in the algorithmic sense. It will also be employed for testing the presence of and estimating the degree of fine-tuning (FT) of a search algorithm or agent who brings about $X$. To this end, we will 
introduce a specificity function $f$, which quantifies, in terms of $f(x)$, how specified an outcome $x\in\Omega$ is. The target $A$, on the other hand, is the set of highly specified states, that is, all states with a degree of specificity that exceeds a given threshold. Then $I^+$ in \re{I+} is a test statistic for testing whether an algorithm has a much larger probability of reaching the set of highly specified states or not, compared to a random search. This is a test of FT, since reaching the target corresponds to specificity (a), whereas reaching it with much higher probability than expected by chance corresponds to (b).  

To calculate $I^+$, the distributions $P_0$ and $P$ of the random search algorithm under $H_0$ and $H_1$, respectively, need to be defined. As mentioned above, the null distribution $P_0$ is typically chosen according to some criterion, such as a maximizer of entropy, possibly with some extra constraints on moments for unbounded $\Omega$, which was the strategy implemented in \cite{DiazHossjerMarks2021, DiazHossjerMarks2022}. Another possibility is to choose $P_0$ as the equilibrium distribution of a Markov chain that models the dynamics of the system under the null hypothesis, for instance an evolutionary process with no external input. In general $P_0=P_{0\xi t}$ involves a number of nuisance parameters $\xi$, and sometimes also the time point $t$ when an algorithm, that does not make use of external information, stops.  The choice of $P=P_{\theta \xi t}$ is problem specific, and it possibly involves the nuisance parameters $\xi$ of the null distribution, the time point $t$ when the algorithm stops, as well as tuning parameters $\theta$ that correspond to infusing background know\-ledge into the search problem. Therefore, in its most general form, the actinfo \re{I+} is a function $I^+=I^+(\theta,\xi,t)$ of the tuning parameters $\theta$, the nuisance parameters $\xi$, and the time point $t$. 

This general framework has many applications, based on different choices of $f$, $A$, $P_0$, and $P$. For some models, $f$ is a binary function that quantifies functionality, so that $A$ is the set of objects of a certain type (e.g., proteins, protein complexes or cellular networks) that are functional among the set $\Omega$ of all such objects. Another possibility is to choose $A$ as the set of populations whose (expected) fitness exceeds a given threshold. In this setting $P_{0\xi t}(A)$ corresponds to the probability that a randomly chosen object or population would reach target $A$ of high fitness at time $t$, given that no background knowledge of the specificity function $f$ is used to generate $X$. The functional information $I_{f_0}=-\log P_{0\xi t}(A)$ corresponds to the amount of external information that an algorithm infuses, given that it brings about $X$ so that $A$ happens with certainty ($P(A)=1$) within time $t$. In this case the object or population is finely tuned when $I_{f_0}$ is large enough. More generally, we say that an evolutionary algorithm that generates $X\sim P=P_{\theta \xi t}$ after $t$ time steps is finely tuned when $I^+(\theta,\xi,t)$ is large enough. Typically, $\theta$ involves selection parameters that determine to which extent a population evolves towards higher fitness.    

The unified treatment of search problems and FT of this paper, is organized as follows: Section \ref{Sec:AI} introduces the specification function $f$ and the set $A$ of highly specified states. Section \ref{Sec:ActiinfoTilt} introduces a class of probability distributions $P=P_\theta$ for which the specificity function $f$ is used to exponentially tilt the null distribution $P_0$, so that outcomes with high specificity are more likely to occur, and with a scalar tuning parameter $\theta$ of $P_\theta$ that corresponds to the amount of exponential tilting. A proof is presented that it is possible to obtain a Metropolis-Hastings type Markov chain in discrete time $t=0,1,2,\ldots$, whose outcome $X=X_t$ at time $t$ has the aforementioned exponentially tilted distribu\-tion under equilibrium, that is, when $t$ is large. The corresponding actinfo $I^+(\theta,t)$ is shown to increase monotonically with $t$ towards an equilibrium limit. The actinfo of a search algorithm $X=X_{t\wedge T}$ that stops at time $T$, when the  targeted set $A$ of highly specified states has been reached, is also shown to increase more rapidly. Section \ref{Sec:Test} introduces various nonparametric and parametric estimators of actinfo, and corresponding tests of FT, when $n$ repeated and independent outputs of the search algorithm are available. In particular, large deviations theory is used to prove that the significance levels of these tests, i.e.\ the probability to detect FT under $H_0$, goes to zero at an exponential rate when the sample size $n$ increases. Section \ref{Sec:Ex} presents a number of examples from cosmology, student learning, reinforcement learning, and population genetics, that illustrate our approach. A discussion in Section \ref{Sec:Disc} follows, whereas proofs and further details about the models, are presented in Section \ref{Sec:Proofs}.    

\section{Specificity and target}\label{Sec:AI}


Consider a function $f:\Omega\to\mathbb R$, and assume that the objective of the search algorithm, or the agent that brings about $X$, is to find regions in $\Omega$ where $f$ is large. The rationale for this is an independent {\it specification}, where a more specified state $x\in\Omega$ corresponds to a larger $f(x)$. It is further assumed that the target set in \eqref{I+} has the form
\begin{align}\label{A}
	A = \{x\in\Omega; \, f(x)\ge f(x_0)\}.
\end{align}

This implies that the purpose of the search algorithm or the agent is to bring about an $X$ that is highly specified. We will refer to $f$ as a specificity function of the agent or an objective function of the search algorithm. For instance, in cosmological FT, $x$ is the value of a particular constant of nature and the specificity function equals
\begin{equation}
f(x) = 1_{\{x\in A\}},
\lb{fbinary}
\end{equation}
where $1_{\{\cdot\}}$ is the indicator function. That is, $f$ has a binary range, with $f(x)=1$ and 0 corresponding to whether $x$ permits a universe with life or not. From this, $A$ is the LPI of this constant if $f(x_0)=1$. Moreover, $X$ is the value of this constant of nature for a randomly generated universe, with a distribution that either incorporates external information ($H_1$) or not ($H_0$). 

In the context of proteins, $x$ is taken to be an amino acid sequence, whereas $f(x)$ in \re{fbinary} quantifies whether the protein that the amino acid corresponds to is functional (1) or not (0). For instance, $X$ could be the outcome of a random evolutionary process, the goal of which is to generate a functioning protein, and this process either makes use of external information ($H_1$) or not ($H_0$). Several other applications are given in Section \ref{Sec:Ex}, including a more refined biological example, where $x$ corresponds to a protein complex or a molecular machine.     

\subsection{Interpretation of target}

There are at least two ways of interpreting $x_0$, and hence also the target set $A$. According to the first interpretation, $x_0$ is the outcome of random variable $X^\prime\in\Omega$; that is, the outcome of a first search. Suppose $X$ is another random variable that represents a second (possibly future) search, independent of $X^\prime$. Then, if we condition on the outcome $x_0$ of the first search,     
the actinfo $I^+$ in (\ref{I+}) is the log likelihood ratio for the event that the second search variable $X$ is {\it at least as specified} as the observed value $f(x_0)$ of the first search. 

There is no need to associate $x_0$ in \eqref{A} with a first search variable $X^\prime$ though. Instead, some {\it apriori information} may be used to define which values of $f$ represent a high amount of specificity. This gives rise to the second interpretation of $x_0$, according to which $x_0$ is used for defining outcomes with a high and low degree of specificity, using $f_0=f(x_0)$ as a cutoff. According to this interpretation, the two sets $A$ in \eqref{A} and its complement
\begin{align*}
	A^c=\Omega\setminus A = \{x; \, f(x) < f(x_0)\}
\end{align*}
represent a dichotomization of specificity, so that $A$ and $A^c$ consist of all states with high and low specificity respectively. With this interpretation of $x$, $I^+$ is the log likelihood ratio for testing FT based on the search variable $X$. In particular, suppose that the specificity function $f$ is bounded, i.e.
\begin{align}\label{fmax}
	f_\text{max} = \max_{x\in\Omega} f(x) < \infty.
\end{align}
Then the most stringent definition of high specificity, 
\begin{align}\label{f0fmax}
	f_0 = f_\text{max},
\end{align}
only regards outcomes with a maximal value of $f$ as highly specified, so that
\begin{align}\label{Omegamax} 
	A = \Omega_\text{max} = \{x\in\Omega; \, f(x)=f_\text{max}\}.
\end{align}
Note that (\ref{fbinary}) is a special case of (\ref{Omegamax}).

\section{Active information for exponentially tilted sys\-tems}\label{Sec:ActiinfoTilt}

Throughout Section \ref{Sec:ActiinfoTilt}, $\xi$ is assumed to be known and the null distribution does not involve any time index $t$. Therefore, $P_0$ is known, whereas $P=P_{\theta t}$ involves the tuning parameters $\theta$ and the time index $t$. It will further be assumed in Sections \ref{Sec:ExpTilt}-\ref{Sec:MH} that the system is in equilibrium, so that the time index $t$ can be dropped also under $H_1$ ($P=P_\theta$). 

\subsection{Exponential tilting}\lb{Sec:ExpTilt}

Let $P_\theta$ be an exponentially tilted version of $P_0$ for some scalar tuning parameter $\theta>0$, which will also be called tilting parameter. Exponential tilting is often used for rare events simulation \cite{AsmussenGlynn2007, Siegmund1976}. Here $f$ is used to define the tilted version of $P_0$ as  
\begin{align}\label{Ptheta}
	P_\theta(x) = \frac{e^{\theta f(x)}}{M(\theta)} P_0(x),
\end{align}
with 
\begin{align}\label{Mtheta}
	M(\theta) = \sum_{x\in\Omega} e^{\theta f(x)} P_0(x)
\end{align}
a normalizing constant assuring that $P_\theta$ is a probability measure. For finite sample spaces $\Omega$, we interpret $P_0(x)$ and $P_\theta(x)$ as probability masses, whereas for continuous sample spaces they are probability densities, and the sum in \eqref{Mtheta} is replaced by an integral. The larger the tilting parameter $\theta>0$, the more the probability mass of $P_\theta$ concentrates on regions of large $f$. In particular, $P_\infty$, the weak limit of $P_\theta$ as $\theta\to\infty$, is supported on (\ref{Omegamax}) whenever \eqref{fmax} holds. 

The parametric family 
\begin{align}\label{cP}
	\mathcal P = \{P_\theta; \, \theta \ge 0\}
\end{align}
of distributions is an exponential family \cite[Section~1.5]{LehmannCasella1998}, and each $P_\theta\in \mathcal P$ gives rise to a separate version of actinfo. This is summarized in  the following proposition:

\begin{prop}\label{Prop:theta}
Suppose the target set $A$ is defined as in \eqref{A} for some $x_0\in\Omega$ such that $P_0(A)>0$. Then $P_\theta(A)$ is a strictly increasing function of $\theta\ge 0$ with $P_\infty(A)=1$. Consequently, the actinfo 
\begin{align}\label{I+theta} 
	I^+(\theta) = \log \frac{P_\theta(A)}{P_0(A)}
\end{align}
is a strictly increasing function of $\theta\ge 0$, with $I^+(0)=0$ and $I^+(\infty) = I_{f0}=-\log P_0(A)$.  
\end{prop} 
  
The intuitive interpretation of Proposition \ref{Prop:theta} is that the larger $\theta$ is, the more problem specific knowledge is infused into $P_\theta$ in terms of shifting probability mass towards regions in $\Omega$ where $f$, the specificity function, is large. 

\subsection{Metropolis-Hastings systems with exponential tilting equilibrium}\lb{Sec:MH}

Inspired by Markov Chain Monte Carlo methods \cite{RobertCasella2010}, consider a Markov chain $X_0,X_1,\ldots \in\Omega$ for which $P_\theta$ is the equilibrium distribution. Consequently, if $P = P_\theta$ (that is, under the alternative hypothesis $H_1$ in (2) when $\theta>0$), $X=X_t$ may be interpreted as the outcome of an algorithm after $t$ iterations, provided $t$ is so large that equilibrium has been reached. The assumption is made that this algorithm knows $f$ and tries to explore the whole state space $\Omega$. If the Markov chain has an equilibrium distribution \eqref{Ptheta}, this corresponds to an algorithm that favors jumps towards regions of large $f$ when $\theta>0$, more so the higher the value of $\theta$ is. In more detail, the transition kernel of the chain is an instance of the well-known Metropolis-Hastings (MH) algorithm \cite{Hastings1970,MetropolisEtAl1953}, which is closely related to simulated annealing \cite{KirkpatrickEtAl1983}. This kernel has a probability or density
\begin{align}\label{pi}
	\pi_\theta(x,y) = r_\theta(x)\delta(x,y) + \alpha_\theta(x,y)q(x,y)
\end{align}
for jumps from $x$ to $y$, where $\delta(x,\cdot)$ is a point mass at $x\in\Omega$, $q(x,\cdot)$ is a proposal distribution of jumps from a current position $x$ of the Markov chain, 
\begin{align}\label{alpha}
	\alpha_\theta(x,y) = \min\left[ 1, \frac{e^{\theta f(y)}P_0(y)q(y,x)}{e^{\theta f(x)} P_0(x)q(x,y)} \right]
\end{align} 
is the probability of accepting a proposed move from $x$ to $y$, whereas 
\begin{align}\label{rtheta}
	r_\theta(x) = 1 - \sum_{y\in\Omega} \alpha_\theta(x,y) q(x,y) 
\end{align}
is the probability that the Markov chain rejects a proposed move away from $x$ (for continuous sample spaces $q(x,\cdot)$ is a probability density and then the sum in \eqref{rtheta} is replaced by an integral). The transition of the Markov chain from $X_t=x$ to the next state $X_{t+1}$ is described in two steps as follows. First a candidate $Y\sim q(x,\cdot)$ is proposed. Then in the second step this candidate is either accepted with probability $\alpha_\theta(x,Y)$, so that $X_{t+1}=Y$, or it is rejected with probability $1-\alpha_\theta(x,Y)$, so that $X_{t+1}=X_t$. It is well known that $P_\theta$ is the equilibrium distribution of this Markov chain whenever it is irreducible; that is, provided the proposal distribution $q$ is defined in such a way that moving between any pair of states in $\Omega$ in a finite number of steps is possible \cite[pp.~243-245]{Ross2003}. 

In particular, if $q$ is symmetric and $P_0$ is uniform, then a proposed upward move with $f(Y)>f(x)$ and $P_\theta(Y)>P_\theta(x)$ is always accepted, whereas a proposed downward move with $f(Y)<f(x)$ is accepted with probability $P_\theta(Y)/P_\theta(x)$. The Markov chain only makes local jumps if $q(x,\cdot)$ puts all its probability mass in a small neighborhood of $x$, for any $x\in\Omega$. At the other extreme is a chain with the global proposal distribution $q(x,\cdot)\sim P_\theta$ for any $x\in\Omega$; all proposed jumps of this chain are then accepted ($\alpha(x,y)=1$), and $\{X_t\}_{t=1}^\infty$ is a sequence of independent and identically distributed (i.i.d.) random variables with $X_t\sim P_\theta$.

\subsection{Active information for Metropolis-Hastings sys\-tems in non-equilibrium}\label{Sec:Noneq}

Suppose for simplicity that the sample space $\Omega$ is finite, and that the states in $\Omega$ are listed in some order. Let 
\begin{align}\label{bP0}
	\textbf P_0 = (P_0(x); \, x\in\Omega)
\end{align}
be a row vector of length $|\Omega|$ with all the null distribution probabilities, and let 
\begin{align}\label{bPi}
	\mathbf \Pi_\theta = \left(\pi_\theta(x,y); x,y\in\Omega\right)
\end{align}
be a square matrix of order $|\Omega|$ that defines the transition kernel of the Markov chain $\{X_t\}_{t=0}^\infty$ of Section \ref{Sec:MH}. If $X_0\sim P_0$, then by the Kolmogorov-Chapman equation $X_t\sim P_{\theta t}$, where 
\begin{align}\label{bPiConv}
	(P_{\theta t}(x); x\in\Omega) = \mathbf P_{\theta t} = \mathbf P_0\mathbf \Pi_\theta^t.
\end{align}
Hence, if $P=P_{\theta t}$, then $X=X_t$ corresponds to observing the Markov chain at time $t$, under the alternative hypothesis $H_1$ in \eqref{H0H1}. Some basic properties of the corresponding actinfo are summarized in the following proposition: 

\begin{prop}\label{Prop:thetat}
Suppose $X=X_t$ is obtained by iterating $t$ times a Markov chain with initial distribution \eqref{bP0} and transition kernel \eqref{bPi}. The actinfo then equals  
\begin{equation}
I^+(\theta,t) = \log \frac{P_{\theta t}(A)}{P_0(A)} = \log \frac{\mathbf P_0\mathbf \Pi_\theta^t\mathbf v}{\mathbf P_0\mathbf v},
\label{I+thetat}
\end{equation}
where $\mathbf v$ is a column vector of length $|\Omega |$ with ones in positions $x\in A$ and zeros in positions $x\in A^c$. In particular, $I^+(\theta,0)=0$ and 
\begin{equation}
\lim_{t\to\infty} I^+(\theta,t) = I^+(\theta).
\label{I+Conv}
\end{equation}
\end{prop}

Therefore, $I^+(\theta,t)>0$ corresponds to knowledge of $f$ being used to generate $t$ jumps of the Markov chain, under the alternative hypothesis $H_1$ in \eqref{H0H1}.  

\subsection{Active information for Metropolis-Hastings sys\-tems with stopping}\label{Sec:Stop}

In Section \ref{Sec:Noneq}, $P\sim P_{\theta t}$ was obtained by starting a random search with null distribution $P_0$, and then iterating the Markov chain of Section \ref{Sec:MH} $t$ times. However, knowledge of $f$ can be utilized even more and stop the Markov chain if the target $A$ in \eqref{A} is reached before time $t$. This can be formalized by introducing the stopping time 
\begin{align}\label{T}
	T = \min\{t\ge 0; \, X_t\in A\}
\end{align}
and letting
\begin{align}\label{Pthetats}
	P_{\theta t s}(x) = P(X_{t\wedge T} =x)
\end{align}
be the probability distribution of the stopped Markov chain $X_{t\wedge T}$, with the last index $s$ in \eqref{Pthetats} being an acronym for stopping. In particular,
\begin{align}\label{Ptheta ts2}
	P_{\theta ts}(A) = \sum_{x\in A} P_{\theta t s}(x) = P(T\le t)
\end{align}
is the probability of reaching the target $A$ for the first time after $t$ iterations or earlier. The theory of phase-type distributions can then be used to compute the target probability $P_{\theta ts}(A)$ in \eqref{Pthetats} \cite{AsmussenEtAl1996,Neuts1981}. To this end, clump all states $x\in A$ into one absorbing state, and decompose the transition kernel  in \eqref{bPi} according to
\begin{align}\label{bPi2}
	\mathbf \Pi_\theta = 
	\left(
	\begin{array}{cc}
		\mathbf \Pi_\theta^\text{na} & \mathbf \Pi^\text{na,a}_\theta \\
		\mathbf 0 & 1 
	\end{array}
	\right),
\end{align}
where $\mathbf \Pi^\text{na}_\theta$ is a square matrix of order $|A^c|$ containing the transition probabilities between all non-absorbing states in $A^c$, whereas $\mathbf \Pi^\text{na,a}_\theta$ is a column vector of length $|A^c|$ with transition probabilities $\pi(x,A)$ from all the non-absorbing states $x\in A^c$ into the absorbing state $A$. Moreover, $\mathbf P_0^\text{na}=\left(P_0(x); \, x\in A^c\right)$ is a row vector of length $|A^c|$ that is the restriction of the start-distribution $\mathbf P_0$ in \eqref{bP0} to all non-absorbing states. Then
\begin{align}\label{Pthetats2}
	P_{\theta ts}(A) = 1 - \mathbf P_0^{\text{na}} (\mathbf \Pi^{\text{na}}_\theta)^{t}\mathbf 1,
\end{align}
where $\mathbf 1$ is a column vector of $|A^c|$ ones.

The actinfo $I^+_s$ of a search procedure with stopping is thus defined: 

\begin{prop}\label{Prop:thetats}
Suppose $X=X_t$ is obtained by iterating a Markov chain with initial distribution \eqref{bP0} and transition kernel \eqref{bPi} (for some $\theta\ge 0$) at most $t$ times, and stopping whenever the set $A$ is reached. Then the actinfo is given by
\begin{align}\label{I+thetats}
	I^+_s(\theta,t) = \log \frac{P_{\theta ts}(A)}{P_0(A)} = \log \frac{1 - \mathbf P_0^{\text{na}} (\mathbf \Pi^{\text{na}}_\theta)^{t}\mathbf 1}{\mathbf P_0\mathbf v},
\end{align}
with $\mathbf P_0$ and $\mathbf v$ as in Proposition \ref{Prop:thetat}, whereas $\mathbf P_0^{\text{na}}$, $\mathbf \Pi_\theta^{\text{na}}$, and $\mathbf 1$ are defined below \eqref{bPi2} and \eqref{Pthetats2}. This actinfo satisfies 
\begin{align}\label{I+Ineq}
	I^+_s(\theta,t) \ge I^+(\theta,t)
\end{align}
and  $I_s^+(\theta,t)$ is a non-decreasing function of $t$ such that 
\begin{align}\label{I+Lim}
	\lim_{t\to\infty} I_s^+(\theta,t) = I_{f0}
\end{align}
and 
\begin{align}\label{I+ET}
	\sum_{t=0}^\infty \left( 1 - P_0(A)e^{I^+_s(\theta,t)} \right) = E(T).
\end{align}
\end{prop}

Inequality \eqref{I+Ineq} states that, for a search procedure with $t$ iterations, knowledge about $f$ that is used for {\it stopping} the Markov chain in \eqref{bPi} will increase the actinfo, regardless of whether knowledge about $f$ was used ($\theta >0$) or not ($\theta=0$) when {\it iterating} the Markov chain. Equation \eqref{I+Lim} is a consequence of the fact that the target $A$ is reached eventually with probability 1, so that the actinfo of a search procedure with stopping equals the functional information $I_{f0}=-\log P_0(A)$ after many iterations of the Markov chain. Moreover, equation \eqref{I+ET} tells that the rate at which $P_0(A)e^{I^+_s(\theta,t)}$ approaches 1 is determined by the expected waiting time $E(T)$ of reaching the target.  

From Proposition \ref{Prop:thetats}, actinfo for a system with stopping is closely related to the phase-type distribution of the waiting time $T$ until the target is reached. This has been studied in \cite{HossjerBechlyGauger2021}, in the context of gene expression of a number of genes, with $x$ the collection of regulatory regions of all these genes.

\section{Estimating active information and testing fine-tuning}\label{Sec:Test}

Suppose the random search algorithm is repeated independently, under the same conditions, $n$ times. For instance, suppose $\{X_{it}\}_{t=0}^\infty$ correspond to independent realizations $i=1,\ldots,n$ of a search algorithm. The outomes of these independent searches are either $X_i=X_{it}$ or $X_i = X_{i,t\wedge T_i}$, for $i=1,\ldots,n$, depending on whether the search algorithm is stopped at a fixed time point $t$ or at random time points $\{T_i\}_{i=1}^n$. In either case, an output of i.i.d random variables 
\begin{equation}
X_1,\ldots,X_n\sim Q
\label{XiQ}
\end{equation}
is obtained. These repeated outcomes of the search algorithm will be used to test for and estimate the degree of fine-tuning. The methodology depends on whether the null distribution $P_0$ is known or involves unknown nuisance parameters.

\subsection{Null distribution known}\label{Sec:H0known}

Suppose the null distribution $P_0$ is known. The sample in \re{XiQ} is then used for testing between the two hypotheses 
\begin{align}\label{H0H12}
	\begin{array}{ll}
		H_0: & Q = P_0,\\
		H_1: & Q \in {\mathcal P}_1,
	\end{array}
\end{align}
with
\begin{equation}
{\mathcal P}_1 = \{P; \, P(A) \ge p_{\text{min}}\}
\label{CP1}
\end{equation}
the set of distributions that correspond to fine-tuning. 
Suppose an estimate $\hat Q(A)$ of the probability that $X\in A$ is computed from data (\ref{XiQ}), with an associated empirical actinfo
\begin{align}\label{hI+}
	\hat I^+ = \hat I^+_n = \log \frac{\hat Q(A)}{P_0(A)}.
\end{align}

If $\hat Q(A)$ is a consistent estimator of $Q(A)$, then for large sample sizes $\hat I^+$ will be close to
\begin{align}\label{I+Q}
	I_Q^+ = \log \frac{Q(A)}{P_0(A)},
\end{align}
which equals 0 under $H_0$ and $I^+=I^+_P$ under $H_1$, for some particular $P\in {\mathcal P}_1$. To test $H_0$ against $H_1$, 
\begin{align}\label{TestAI}
	\text{reject } H_0 \text{ when } \hat I^+\ge I_\text{min},
\end{align}
where $I_{\text{min}}$ is a pre-specified lower bound on the range of values of the actinfo that correspond to FT. 

\subsubsection{Nonparametric estimator and test}

From Section \ref{Sec:ActiinfoTilt}, $P=P_\theta$, $P=P_{\theta t}$ or $P=P_{\theta ts}$ involves the tilting parameter $\theta$, and possibly also the number of iterations $t$ of the algorithm and a stopping time $T$.  In this section, no other assumption than $P\in {\mathcal P}_1$ is made on $P$, and a nonparametric version of the empirical actinfo is used. The fraction 
\begin{align}\label{hQA}
	\hat Q(A) = \frac{1}{n}\sum_{i=1}^n 1_{\{X_i\in A\}}  
\end{align}
of random searches that fall into $A$ is used as an estimate of $Q(A)$. Therefore, (\ref{hQA}) only requires knowledge of the set $A$, not of the function $f$. 

The following result establishes asymptotic normality of the nonparametric version of the estima\-tor $\hat I^+$ in \re{hI+}. Moreover, large deviations \cite{Varadhan1984} are used to show that the significance level of the nonparametric version of the FT test \re{TestAI} goes to zero exponentially fast with $n$: 

\begin{prop}\label{Prop:LD}
	Suppose the empirical actinfo $\hat I_n^+$ in \eqref{hI+} is computed non-parametrically, using \eqref{hQA} as an estimate of the target probability $Q(A)$. Then $\hat I_n^+$ is an asymptotically normal estimator of $I^+_Q$ in \eqref{I+Q}, in the sense that
	\begin{align}\label{AsNNP}
		\sqrt{n}(\hat I^+_n-I^+_Q)  \stackrel{\cal L} \longrightarrow N(0,V)\text{ as } n\to\infty,
	\end{align}
	where $ \stackrel{\cal L} \longrightarrow$ refers to convergence in distribution, and 
	\begin{align}\label{VNP}
		V = \frac{1-Q(A)}{Q(A)}
	\end{align}
	is the variance of the limiting normal distribution. The significance level of the test \eqref{TestAI} for fine-tuning, with threshold $I_{\text{min}}$, satisfies
	\begin{align}\label{LD}
		\lim_{n\to\infty} - \frac{\log \left(P_{H_0}(\hat I_n^+ \ge I_\text{min}) \right)}{n} = C,
	\end{align}
	where 
	\begin{align}\label{C}
		C =  p_\text{min} \log \frac{p_\text{min}}{P_0(A)} + (1-p_\text{min})\log \frac{1-p_\text{min}}{1-P_0(A)}
	\end{align}
	is the Kullback-Leibler divergence between Bernoulli distributions with success probabilities $p_\text{min}=P_0(A)\exp(I_{\text{min}})$ and $P_0(A)$ respectively. 
\end{prop}

\begin{rem}
	The conclusion of Proposition \ref{Prop:LD} is that the probability of observing actinfo that corresponds to fine-tuning by chance decays at rate $e^{-Cn}$ when the sample size $n$ gets large. 
\end{rem}

\subsubsection{Parametric estimator and test}

Suppose there is a priori knowledge that $P$ is close to the parametric exponential family $\mathcal P$ of distributions in \eqref{Ptheta}-\eqref{cP} for some value $\theta>0$ of the tilting parameter. A parametric test of actinfo is naturally defined. For this, compute first the maximum likelihood estimate 
\begin{align}\label{hth}
	\hat \theta = \hat \theta_n = \arg\max_{\theta\ge 0} \sum_{i=1}^n \log P_\theta (X_i)
\end{align}
of $\theta$, and use it to define a parametric estimate 
\begin{align}\label{hQA2}
	\hat Q(A) = P_{\hat \theta}(A)
\end{align}
of the target probability $Q(A)$ that is inserted into \eqref{hI+} to define a parametric version of the empirical actinfo $\hat I^+$. As opposed to (\ref{hQA}), the estimate (\ref{hQA2}) requires full knowledge of $f$. 

To analyze the properties of the estimator \eqref{hI+} and test \eqref{TestAI}, introduce 
\begin{align}\label{thetaast}
	\theta^\ast = \arg\min_{\theta\ge 0} D_{KL}(Q~\Vert~P_\theta),
\end{align}
where 
\begin{align}\label{DKL}
	D_{KL}(Q~\Vert~P_\theta) = \sum_{x\in\Omega} Q(x) \log \frac{Q(x)}{P_\theta(x)}
\end{align}
is the Kullback-Leibler divergence between $Q$ and $P_\theta$. From \eqref{thetaast}, $P_{\theta^\ast}$ is the distribution in $\mathcal P$ that best approximates $Q$. In particular, $\theta^\ast=\theta$ if $Q\in\mathcal P$ and $Q=P_\theta$ for some $\theta\ge 0$. 

The following proposition shows that $\hat I^+$ is an asymptotically normal estima\-tor of $I^+(\theta^\ast)$ in \eqref{I+theta}, which differs from $I_Q^+$ in \eqref{I+Q} whenever $Q\notin \mathcal P$. Moreover, the proposition also provides large sample properties of the significance level of the test for actinfo:

\begin{prop}\label{Prop:Param}
	Suppose the empirical actinfo $\hat I_n^+$ in \eqref{hI+} is computed parametrically, using an estimate \eqref{hQA2} of the target probability $Q(A)$. Then $\hat I_n^+$ is an asymptotically normal estimator of $I^+(\theta^\ast)$, in the sense that
	\begin{align}\label{AsNP}
		\sqrt n \left(\hat I^+_n-I^+(\theta^\ast) \right) \stackrel{\cal L} \longrightarrow N(0,V)\text{ as }n\to\infty,
	\end{align}
	where the variance of the limiting normal distribution is given by
	\begin{align}\label{VPar}
		V = \frac{\text{Cov}^2_{P_{\theta^\ast}} \left[ f(X)I(f(X)\ge f_0) \right] \text{Var}_Q \left[ f(X) \right]}{P_{\theta^\ast}^2(A) \text{Var}^{\,\, 2}_{P_{\theta^\ast}}			[f(X)]}.
	\end{align}
	Moreover, the significance level of the parametric test for fine-tuning, based on \eqref{TestAI} and \eqref{hQA2}, satisfies  
	\begin{align}\label{LD2}
		\lim_{n\to\infty} - \frac{\log \left[ P_{H_0} \left(\hat I_n^+ \ge I_\text{min} \right) \right]}{n} = C,
	\end{align}
	for  
	\begin{align}\label{C2}
		C = \sup_{\phi>0} \left\{\phi E_{P_\text{min}}[f(X)] - \log M(\phi)\right\},
	\end{align}
	where $P_\text{min} = P_{\theta_\text{min}}$, $\theta_\text{min}<\theta^\ast$ is the solution of $P_{\theta_\text{min}}(A)=p_\text{min}=P_0(A)\exp(I_{\text{min}})$, $M(\phi)$ is given by \eqref{Mtheta}, whereas $p_{\text{min}}$ is defined in \re{CP1}.  
\end{prop}

\subsubsection{Comparison between nonparametric and parametric estimates of actinfo}

The two versions of empirical actinfo are complementary. The nonparametric version is preferable in the sense that it makes less assumptions about the distribution $P$ of the random algorithm under $H_1$, and in particular it is a consistent estimator of $I_Q^+$ in \eqref{I+Q}. The parametric version of $\hat I^+$, on the other hand, is preferable when $nQ(A)$ is small, since it makes use of all data in order to estimate $Q(A)$, although it is not a consistent estimator of $I_Q^+$ when $Q\notin \mathcal P$. The asymptotic variances in \eqref{VNP} and \eqref{VPar}, as well as the rates of exponential significance level decrease in \eqref{C} and \eqref{C2}, agree when $Q=P_{\theta^\ast}$ and $f(x) = f_0 1_{\{x\in A\}}$, which is a special case of \eqref{f0fmax}.     

\subsection{Null distribution unknown}\label{Sec:H0unknown}

Suppose the null distribution $P_0=P_{0\xi}$ involves an unknown nuisance parameter $\xi\in\Xi$. The objective is then to test the two hypotheses
\begin{align}\label{H0H13}
	\begin{array}{ll}
		H_0: & Q \in {\mathcal P}_0,\\
		H_1: & Q \in {\mathcal P}_1,
	\end{array}
\end{align}
where the set of distribution under the null and alternative hypotheses equals 
\begin{equation}
{\mathcal P}_0 = \{P_{0\xi}; \, \xi\in\Xi\}
\end{equation}
and \re{CP1} respectively. 

\subsubsection{One sample available}

The actinfo
\begin{equation}
I_Q^+ = I_Q^+(\xi) = \log \frac{Q(A)}{P_{0\xi}(A)}
\lb{I+nuisance}
\end{equation}
cannot be consistently estimated if only one sample \re{XiQ} is available. The best that can be done is to estimate a lower bound 
\begin{equation}
\hI^+ = \hI_n^+ = \log \frac{\hQ(A)}{P_{0\text{max}}(A)}
\lb{hI+max}
\end{equation}
of $I^+$, with $P_{0\text{max}}(A)$ defined in \re{P0max} and $\hQ(A)$ an estimate of $Q(A)$. This estimator will have an asymptotic bias 
\begin{equation}
B = I^+_Q(\xi^\ast) - I_Q^+ = \log \frac{P_{0\xi}(A)}{P_{0\text{max}}(A)} \le 0, 
\lb{bias}   
\end{equation}
with $\xi^\ast$ the nuisance parameter that maxizes $P_{0\xi}(A)$ \cite{HossjerEtAl2022}. For the numerator of \re{hI+max} either the nonparametric estimate of $Q(A)$ in \re{hQA} can be used, or a parametric class
$$
{\mathcal P} = \{P_{\theta\xi}; \, \theta\in\Theta,\xi\in\Xi\}
$$
of distributions can be used that involves a tuning parameter vector $\theta$ and a vector of nuisance parameters $\xi$.  
If $Q$ is thought to be close to $\cP$, the parametric estimate
\begin{equation}
\hQ(A) = P_{\hat{\theta}\hat{\xi}}(A)
\lb{hQA3}
\end{equation}
of $Q(A)$ is used that generalizes \re{hQA2}, with 
\begin{equation}
(\hat{\theta},\hat{\xi}) = \arg\max_{\theta,\xi} \sum_{i=1}^n \log P_{\theta\xi}(X_i).
\lb{hthhxi}
\end{equation}
When the sample size $n$ tends to infinity, the estimator \re{hthhxi} will converge to
\begin{equation}
(\theta^\ast,\xi^\ast) = \arg\min_{\theta,\xi} D_{KL}(Q~\Vert~P_{\theta\xi}).
\lb{thetaxiast}
\end{equation}
The following result is an extension of Propositions \ref{Prop:LD}-\ref{Prop:Param}, when nuisance parameters $\xi$ are added and a general type of tuning parameter $\theta$ (not necessarily a scalar tilting parameter) is used: 

\begin{prop}\lb{Prop:Nuisance} Suppose the null distribution $P_0=P_{0\xi}$ involves an unknown parameter $\xi$ and the actinfo $I_Q^+$ in \re{I+nuisance} is estimated by $\hI_n^+$ in \re{hI+max}, using an estimator $\hQ(A)$ of the target probability $Q(A)$ that is either nonparametric (\ref{hQA}) or parametric \re{hQA3}. Given these assumptions, $\hI_n^+$ is an asymptotically normal estimator, in the sense that
\begin{equation}
\sqrt{n}(\hI^+_n - I_Q^+ - B) \stackrel{\cal L} \longrightarrow N(0,V)\text{ as }n\to\infty.
\label{AsNnuisance}
\end{equation}
The asymptotic bias $B$ in \re{AsNnuisance} is defined in \re{bias} whereas the asymptotic variance $V$ is defined in \re{VNP} for the nonparametric estimator of $I_Q^+$, whereas 
\begin{equation}
\begin{array}{rcl}
V &=& E[\psi_{\theta^\ast \xi^\ast}(X)|X\in A] E[\psi_{\theta^\ast \xi^\ast}^\prime (X)]^{-1} E[\psi^T_{\theta^\ast \xi^\ast}(X)\psi_{\theta^\ast \xi^\ast}(X)]\\
&\cdot & E[(\psi_{\theta^\ast \xi^\ast}^\prime)^T(X)]^{-1}E[\psi_{\theta^\ast \xi^\ast}(X)|X\in A]^T
\end{array}
\lb{VPar2}
\end{equation}
for the parametric estimator of $I_Q^+$, with $\psi_{\theta\xi}(x) = d\log P_{\theta \xi}(x)/d(\theta,\xi)$, $(\theta^\ast,\xi^\ast)$ defined as in \re{thetaxiast}, and $T$ refering to matrix transposition. Moreover, the signifi\-cance level of the test \re{TestAI} of FT, with threshold $I_\text{min}$, satisfies   
\begin{align}\label{LD3}
	\lim_{n\to\infty} - \frac{\log \left[ P_{0\xi} \left(\hat I_n^+ \ge I_\text{min} \right) \right]}{n} = C,
\end{align}
with 
\begin{equation}
C = p_\text{min}e^{-B} \log \frac{p_\text{min}e^{-B}}{P_{0\xi}(A)} + (1-p_\text{min}e^{-B})\log \frac{1-p_\text{min}e^{-B}}{1-P_{0\xi}(A)}
\lb{CNonpar}
\end{equation}
for the nonparametric version of the test, with $p_{\text{min}}=P_{0\xi}(A)\exp(I_{\text{min}})$.  For the parametric versions of the FT-test, and in the special case when $\theta$ is a scalar exponential tilting parameter, $C$ is given by \re{C2}, with $P_{\text{min}}=P_{\theta_{\mbox{\tiny min}}\xi}$, and $\theta_{\text{min}}$ the solution of $P_{\theta_{\mbox{\tiny min}}\xi}(A)=p_{\text{min}}e^{-B}$.      
\end{prop}

\begin{rem} The negative bias term $B$ makes the test of FT in Proposition \ref{Prop:Nuisance} more conservative than the tests in Propositions \ref{Prop:LD}-\ref{Prop:Param}. This can be seen, for instance, by comparing the two large deviation rates $C$ in \re{C} and \re{CNonpar}. The rate in \re{CNonpar} is larger, since $p_{\text{min}}$ is multiplied by a term $e^{-B}$. This corresponds to the fact that to falsely reject $H_0$ in Proposition \ref{Prop:Nuisance} is more difficult.   
\end{rem} 

\subsubsection{Two samples available}

In addition to the first sample \eqref{XiQ}, suppose a second sample 
\begin{equation}
X_{01},\ldots,X_{0n_0}\sim P_{0\xi}
\label{XiQ0}
\end{equation}
of $n_0$ i.i.d.\ observations under the null distribution is available. A consistent estimator 
\begin{equation}
\hI^+ = \hI^+_{nn_0} = \log \frac{\hQ(A)}{P_{0\hat{\xi}}(A)} 
\lb{hITwoSamp}
\end{equation}
of $I_Q^+$ in \re{I+nuisance} is then available, with
\begin{equation}
\hat{\xi} = \arg\max_{\xi} \sum_{i=1}^{n_0} \log P_{0\xi}(X_{0i}).
\lb{hxi}
\end{equation}

The following result provides asymptotic properties of the estimator \re{hITwoSamp} of actinfo, and the corresponding test \re{TestAI} of FT with threshold $I_\text{min}$: 

\begin{prop}\label{Prop:TwoSamples}
Suppose the null distribution $P_0=P_{0\xi}$ involves an unknown nuisance parameter $\xi$, and that the active information     
$I_Q^+$ in \re{I+nuisance} is estimated by $\hI_{nn_0}^+$ in \re{hITwoSamp}, making use of two samples \re{XiQ} and \re{XiQ0}, of sizes $n$ and $n_0$, from $Q$ and $P_{0\xi}$ respectively. Assume further that the estimator $\hQ(A)$ of $Q(A)$ is either nonparametric \re{hQA} or parametric \re{hQA3}. If $n,n_0\to\infty$ in such a way that
\begin{equation}
\frac{n}{n_0} \to \la > 0,
\lb{n0nAsympt}
\end{equation}
then 
\begin{equation}
\sqrt{n}(\hI^+_{nn_0}-I_Q^+) \Lto N(0,V_1 + \la V_2),
\lb{hIAsNTwoSamp}
\end{equation}
where
\begin{equation}
V_2 = E[\psi_\xi(X)|X\in A] E[\psi_\xi^T(X)\psi_\xi(X)]^{-1} E[\psi_\xi(X)|X\in A]^T,
\lb{V2TwoSamp}
\end{equation}
and $\psi_\xi(x) = d \log P_{0\xi}(x) / d\xi$. If the nonparametric estimator of $Q(A)$ is used, then $V_1$ equals $V$ in \re{VNP}, whereas if the parametric estimator $Q(A)$ is used, then $V_1$ equals $V$ in \re{VPar2}. The significance level of the test \re{TestAI} of FT, with threshold $I_{\text{min}}$, satisfies the same type of large deviation result \re{LD3} as in Proposition \ref{Prop:Nuisance}, for the nonparametric and parametric versions of the test (in the latter case assuming that $\theta$ is a scalar tilting parameter), but in the definitions of the nonparametric and parametric large deviation rates $C$, the bias term $B=0$.  
\end{prop}

\section{Examples}\label{Sec:Ex}

\begin{exa}[Cosmology \cite{DiazHossjerMarks2021, DiazHossjerMarks2022}] Suppose there is a positive constant of nature $X\in \Omega = \R^+$, a life-permitting interval $A\subset \Omega$, and a specificity function \re{fbinary} that equals 1 inside $A=(a,b)$ and zero elsewhere. The maximum entropy distribution under a first moment constraint $\xi = E(X)$ is exponential with expected value. Consequently, 
$$
P_{0\xi}(A) =  \frac{1}{\xi}\int_a^b e^{-x/\xi}dx.
$$ 
The null and alternative hypotheses for the fine-tuning test are given in \re{H0H13}, where under $H_1$ the agent brings about a life-permitting value of $X$ with probabi\-lity 1 ($P(A)=1$). Only one universe is observed, with a value $X=X_1$ of the constant. Therefore, there is a sample \re{XiQ} of size $n=1$, whereas no null sample \re{XiQ0} is available. Since $X_1\in A$ is life-permitting, $\hQ(A)=1$. The estimate \re{hI+max} of actinfo then simplifies to  
\begin{equation}
\hI^+ = \log \frac{1}{P_{0\text{max}}(A)} = - \log P_{0\text{max}}(A).
\lb{hI+maxCosm}
\end{equation}
Let $x=(a+b)/2$ be the midpoint of the LPI and suppose that half of its relative size $\epsilon = (b-a)/(2x)$ is small. The probability in \re{hI+maxCosm} is then approximated by  
$$
P_{0\text{max}}(A) \approx (b-a) \max_{\xi >0} \frac{e^{-x/\xi}}{\xi} \approx 2\epsilon e^{-1}.
$$
From \re{hI+maxCosm} the estimated actinfo
$$
\hI^+ \approx 1 - \log(\epsilon) - \log(2)
$$
is a monotone decreasing function of $\epsilon$. 
\end{exa}

\begin{exa}[Evaluation of student test scores \cite{HossjerDiazRao2022}] Suppose a number of students perform a test. Let $x=(z,y)=(z_1,\ldots,z_{d-1},y)\in\R^d$ summarize the chararcte\-ris\-tics of a student with covariates $z$ that are used to predict the outcome $y$ of the test. The specificity function $f(x)=x_d=y$ equals the student's test score, and \re{A} corresponds to the set of students that pass the test, with a minimally allowed score of $f_0$. The population of students follows a $(d-1)$-dimensional multivariate normal distribution $Z\sim N(\bm,\bSi)$, where $\bm=(m_1,\ldots,m_{d-1})$ and $\bSi=(\si_{jk})_{j,k=1}^{d-1}$ are known. The conditional distribution of the response follows a multiple linear regression model 
$$
Y|Z=z \sim N\left(\xi_0 + \sum_{j=1}^{d-1} \xi_jz_j + t(\theta_0+\sum_{j=1}^{d-1} \theta_j z_j),\si^2 \right),
$$
for a student with covariate vector $z$ who prepared for the test for a period of length $t$. The nuiscance parameter vector $\xi=(\xi_0,\ldots,\xi_{d-1},\si^2)$ involves the error variance and the regression parameters for students who did not train for the test, whereas the tuning parameter vector $\theta=(\theta_0,\ldots,\theta_{d-1})$ involves the regression parameters that correspond to the effect of preparing for the test. The unconditional distribution of the response is normal, $Y\sim N(\mu,V)$, with
$$
\begin{array}{rcl}
\mu &=& \mu(\theta,\xi,t) = (\xi_0 + t\theta_0) + \sum_{j=1}^{d-1} (\xi_j+t\theta_j)m_j,\\
V &=& V(\theta,\xi,t) = \si^2 + \sum_{j,k=1}^{d-1} (\xi_j+t\theta_j)(\xi_k+t\theta_k)\si_{jk}.
\end{array}
$$
Therefore, the probability, for a randomly chosen student, that studied for the test for a period of length $t$, to pass is 
\begin{equation}
P(A) = P_{\theta\xi t}(A) = P(Y\ge f_0) = 1 - \Phi\left(\frac{f_0-\mu}{\sqrt{V}}\right),
\lb{PAStudent}
\end{equation}
where $\Phi$ is the cumulative distribution function of a standard normal distribution. The null distribution $P_0=P_{0\xi}$ corresponds to putting $t=0$ in \re{PAStudent}. Thus the actinfo 
\begin{equation}
I^+ = I^+(\theta,\xi,t) = \log \frac{1 -\Phi\left((f_0-\mu(\theta,\xi,t)/\sqrt{V(\theta,\xi,t)}\right)}{ 1 -\Phi\left((f_0-\mu(0,\xi,0)/\sqrt{V(0,\xi,0)}\right)}
\label{I+Student}
\end{equation} 
quantifies how much learning, during a period of length $t$, increases the probability of passing the test. To compute an estimate $\hI^+$ of $I^+$ in \eqref{I+Student}, estimates $\hxi$ and $\hth$ of $\xi$ and $\theta$ are needed. This can be done by collecting two training samples, as in \re{hITwoSamp}. Another option is computing least squares estimates of the nuisance and tuning parameters $(\xi,\theta)$ jointly, without bias, from one single data set $\{(t_i,z_i,y_i)\}_{i=1}^n$, provided that the time periods $t_i$ vary, so that all parameters are identifiable.  
\end{exa}

\begin{exa}[Reinforcement learning (RI) \cite{KaelblingLittmanMoore1996}] Consider an agent whose purpose is to maximize the reward $f(x)$ of a trajectory $x$ that he to some extent will be able to control, for a time period of length $t$. At each time point $u$ there are $m$ possible environments $\cS=\{s_1,\ldots,s_m\}$ and $q$ possible actions $\cA=\{a_1,\ldots,a_q\}$ to take. The state space $\cX=\cA^t\times \cS^{t+1}$ consists of all possible trajectories 
$$
x = (a_0,\ldots,a_{t-1},s_0,\ldots,s_t)
$$
of environments and actions, where $s_u$ is the environment and $a_u$ the action taken at time $u$. A corresponding random trajectory is denoted with capital letters 
$$
X = (A_0,\ldots,A_{t-1},S_0,\ldots,S_t).
$$
If the environment of the system is $S_u=s$ at time $u$, and action $A_u=a$ is taken, the probability of moving to environment $s^\prime$ is $P_a(s,s^\prime)=P(S_{u+1}=s^\prime|S_u=s,A_u=u)$, with an instantaneous reward of $R_a(s,s^\prime)$. If future rewards are discounted by a factor $\ga$, the total reward, over a time horizon of length $t$, is
$$
f(x) = \sum_{u=0}^t R_{a_u}(s_u,s_{u+1})\ga^u.
$$
Let $f_0$ be a lower bound for a trajectory's total discounted reward to be acceptable, so that $A$ in \re{A} is the set of all acceptable trajectories.   
The agent takes action according to some {\sl policy} to make the expected total reward of a trajectory as large as possible. To this end consider stationary policies, where the action $A_u$ taken by the agent at each time point $u$ is only determined by the current environment $s_u$, according to some matrix $\bPi=(\pi(s,a);s\in\cS, a\in\cA)$ of transition probabilities $\pi(s,a)=P(A_u=u|S_u=s)$. For a completely random policy
$$
\pi(s,a) = \xi_a; \quad a=1,\ldots,q,
$$
the action is not influenced by the current environment, and it is completely specified by the vector $\xi=(\xi_1,\ldots,\xi_q)$ of nuisance parameters. Thus $P_0(A)=P_{0\xi t}(f(X)\ge f_0)$ is the probability that an ignorant agent with policy determined by $\xi$, will have an acceptable trajectory. An agent who knows the reward function $R_a$ and the dynamics $P_a$ of the environment, will try to take this knowledge into account to formulate a policy that makes the reward as large as possible. A deterministic policy $\theta: \cS \to \cA$ is a function that to each environment takes a unique action, so that
$$
\pi(s,a) = 1_{\{a=\theta(s)\}}.
$$
Thus $P(A)=P_{\theta t}(f(X)\ge f_0)$ is the probability that an agent with deterministic policy $\theta$ obtains an acceptable trajectory. The active information 
\begin{equation}
I^+ = I^+(\theta,\xi,t) = \log
\frac{P_{\theta}(\sum_{u=0}^t R_{A_u}(S_u,S_{u+1})\ga^u\ge f_0)}
{P_{0\xi}(\sum_{u=0}^t R_{A_u}(S_u,S_{u+1})\ga^u\ge f_0)}
\label{I+RI}
\end{equation}
quantifies, on a logarithmic scale, how much more likely it is for an agent with policy $\theta$ to obtain an acceptable trajectory, compared to an ignorant agent with policy $\xi$. The values $\xi$ and $\theta$ are varied during the exploration phase of RI, but they are assumed to be known during the exploitation phase of RI. Suppose we want to compute the actinfo \eqref{I+RI} during the expoitation phase. Since $P_0(A)$ and $P(A)$ are typically unknown, they have to be estimated by Monte Carlo. To this end, assume we have two samples \eqref{XiQ} and \eqref{XiQ0} of $n$ and $n_0$ trajectories available, from $Q=P_{\theta t}$ and $Q=P_{0\xi t}$ respectively. Then $\hI^+$ in \eqref{hITwoSamp} can be used to estimate the actinfo \eqref{I+RI}. 
\end{exa}

\begin{exa}[Molecular machines and Moran models \cite{ThorvaldsenHossjer2020, HossjerBechlyGauger2021, Montanez2018}]\lb{Ex:Moran}
Suppose $\Omega$ consists of all $2^d$ binary sequences $x=(x_1,\ldots,x_d)$ of length $d$, with a null distribution $P_0(x)$ that will be chosen below. The specificity function $f$ is defined as
\begin{align}\label{fMoran}
	f(x) = \left\{
	\begin{array}{ll}
		a |x|, & x\ne (1,\ldots,1),\\
		1, & x=(1,\ldots,1),
	\end{array}\right.
\end{align}
where $|x|=\sum_{i=1}^d x_i$ and $a \le 1/d$ is a fixed parameter. We regard $x$ as a molecular machine with $d$ parts, with $x_i=1$ or 0 depending on whether part $i$ functions or not. The specificity $f(x)$ quantifies how well the machine works, for instance its ability to regulate activity {\sl in vitro} or {\sl in vivo} in a living cell. It is assumed that $f(x)$ is determined by the number $|x|$ of functioning parts, with a maximal value $f_\text{max} = f(1,\ldots,1) = 1$. Using  \eqref{f0fmax}, the most stringent definition of high specificity, it follows that $A = \{(1,\ldots,1)\}$ only contains one element, a molecular machine for which all parts are in shape. The parameter $a$ is crucial. If $0<a\le 1/d$, it follows that a molecular machine works better the more of the parts that are in shape. On the other hand, if $a <0$, then a molecular machine with some parts in shape, but not all, functions worse the more of the parts that are in shape, since all units must work in order for the whole machine to function, and there is a cost $-a$ associated with carrying each part that is in shape, as long as the whole system does not function. 

Each state $x$ is interpreted as a {\it population} of $N$ subjects, all having the same variant $x$ of the molecular machine. With this interpretation, $X=X_t$ is the outcome of a random evolutionary process where all subjects of the population, at any time point $t$, have the same state. But this state may vary over time when all subjects of population simultaneously experience the same change. The question of interest is whether this process can modify the population so that all its members have a functioning molecular machine. A transition of this process from $x$ is caused by a mutation with distribution $q(x,\cdot)$, where $q(x,x)=0$. Suppose a mutation from $x$ to $y$ is possible, i.e., $q(x,y)>0$. A mutation from $x$ to $y$ first occurs in one individual and then it either (momentarily) dies out with probability $1-\alpha_{\theta}(x,y)$ or it (momentarily) spreads to the whole population (gets fixed) with probability 
\begin{align}\label{alphaN}
	\alpha_{\theta}(x,y) = C \cdot \left(\frac{e^{\theta f(y)}P_0(y)q(y,x)}{e^{\theta f(x)}P_0(x)q(x,y)}\right)^{1/2}, 
\end{align}
where
\begin{align}\label{CN}
	C = \left(\max_{x,y} \frac{e^{\theta f(y)}P_0(y)q(y,x)}{e^{\theta f(x)}P_0(x)q(x,y)}\right)^{-1/2}
\end{align}
is a constant assuring that \eqref{alphaN} never exceeds 1, and the maximum is taken over all $x,y$ such that $x\ne y$ and both of $q(x,y)$ and $q(y,x)$ are positive. The Markov chain with transition probabilities \eqref{pi} and acceptance probability \eqref{alphaN} represents the dynamics of the evolutionary process. 

As shown in Section \ref{Sec:Proofs},
the equilibrium distribution of this Markov chain is given by $P_\theta$ in \eqref{Ptheta}. In particular, Propositions \ref{Prop:thetat}--\ref{Prop:thetats} remain valid when the Markov chain \eqref{pi} with acceptance probabilities \eqref{alphaN} are used, rather than \eqref{alpha}. We will interpret
\begin{align}\label{sx}
s(x) = e^{\theta f(x)/N}
\end{align} 
as the selection coefficient or fitness of individuals with a molecular machine of type $x$, that is, $s(x)$ is proportional to the fertility rate of individuals of type $x$. 

The MH-type Markov chain with acceptance probability \eqref{alphaN}--\eqref{CN} represents an evolutionary process that closely resembles a Moran model with selection \cite{Durrett2008, Moran1958a, Moran1958b}, which is frequently used for describing evolutionary processes (see Section \ref{Sec:Proofs}). The Moran model is a continuous time Markov chain for a population with overlapping generations where individuals die at the same rate, and are replaced by offspring of individuals in the population proportionally to their selection coefficients $s(x)$. New types arise when an offspring of parents of type $x$ mutate with probability $\mu(x)$. If the mutation rate is small ($\mu(x)\ll N^{-1}$ for all $x\in\Omega$), then to a good approximation the whole population will have the same type at any point in time, a so called fixed state assumption. 

Even though the Moran model is specified in continuous time, time can be discretized as $t=0,1,2,\ldots$ by only recording the population when individuals die. If individuals die at rate 1, then the next individual dies at rate $N$, so that time is counted in units of $N^{-1}$ generations. The fixed state assumption is motivated by assuming that newborn offspring with a new mutation either dies out or spreads to the whole population (get fixed in the population) right after birth. In this context, $q$ corresponds to the way in which mutations change the type of the individual, whereas $\alpha_\theta = \alpha_{\theta N}$ is the probability of fixation. If $q(x,y)$ is the conditional probability that an offspring of a type $x$ parent mutates to $y$, given that a mutation occurs, then the proposal kernel of the Moran model is 
\begin{align}\label{qMoran}
	q^\text{Moran}(x,y) = 
		\left\{
		\begin{array}{ll}
		\mu(x) q(x,y), & x\ne y,\\
		1-\mu(x), & x=y.
		\end{array}
		\right.
\end{align}
As shown in Section \ref{Sec:Proofs}, the acceptance (or fixation) probability of the Moran model is 
\begin{align}\label{alphaMoran}
	\alpha^{\text{Moran}}_{\theta N}(x,y) \approx \frac{1}{N} \left(1 + \frac{\theta [f(y)-f(x)]}{2}\right) \approx  \frac{1}{N} \left( \frac{e^{\theta f(y)}}{e^{\theta f(x)}} 			\right)^{1/2}
\end{align}
when $\theta [f(y)-f(x)]$ is small. From \eqref{qMoran}-\eqref{alphaMoran}, the Moran model approximates the Metropolis-Hastings kernel with acceptance probabilities \eqref{alphaN}-\eqref{CN}  with good accuracy when i) $\mu(x)\equiv \mu$, ii) $P_0$ is uniform and iii) the proposal kernel $q$ is symmetric (i.e.\ $q(x,y)=q(y,x)$), although the time scales of the two processes are different. More specifically, if i)-iii) hold, a time-shifted version of the Moran model approximates the MH-type model with acceptance probabilities \eqref{alphaN}-\eqref{CN}, so that each time step of the MH-type Markov chain corresponds to $C/\mu$ generations of a Moran model. However, even under assump\-tions i)-iii) the stationary distribution of the Moran model differs slightly from $P_\theta$. 

The proposal kernel $q(x,y)$ is assumed to be local and satisfying 
\begin{align}\label{q}
	q(x,y) = 
	\left\{
	\begin{array}{ll}
		b / [|x| + b(d-|x|)], & y = x + e_j, x_j=0,\\
		1 / [|x| + b(d-|x|)], & y = x + e_j, x_j=1,\\
		0, & \text{otherwise},
	\end{array}
	\right.
\end{align}
where $e_j=(0,\ldots,0,1,0,\ldots,0)$ is a row vector of length $d$ with a 1 in position $j\in\{1,\ldots,d\}$ and zeros elsewhere, whereas $x+e_j$ refers to component-wise addition modulo 2, corresponding to a switch of component $j$ of $x$. A change of component $j$ from 0 to 1 is caused by a beneficial mutation, whereas a change from 1 to 0 corresponds to a deleterious mutation. Consequently, $b>0$ is the ratio between the rates at which beneficial and deleterious mutations occur. 

The kernel $q$ in \eqref{q} is symmetric only when beneficial and deleterious mutations have the same rate ($b=1$). The more general case of asymmetric $q$ is handled differently by the MH-type algorithm and the Moran model. Whereas the MH-type algorithm elevates the acceptance probability \eqref{alphaN} of seldom-proposed states $y$ (those $y$ for which $q(x,y)$ is small for many $x$), this is not the case for the acceptance probability \eqref{alphaMoran} of the Moran model. To avoid that these states $y$ are reached too often by the MH-type algorithm, the null distribution $P_0$ of no selection has to be chosen so that $P_0(y)$ is small for rarely proposed states (whereas the Moran model needs no such correction). Therefore $P_0$ in \eqref{alphaN} will be chosen as the stationary distribution of a transition kernel \eqref{pi} for which $\theta=0$ and all candidates are accepted ($\alpha_{0}(x,y)=1$). That is, if $\tilde{\mathbf \Pi}_0$ refers to the transition matrix of such a Markov chain, the initial distribution $\mathbf P_0$ in \eqref{bP0} is chosen as the solution of 
\begin{align}\label{P0MHtype}
	\left\{
	\begin{array}{ll}  
		\mathbf P_0 = \mathbf P_0 \tilde{\mathbf \Pi}_0,\\
		\sum_{x\in\Omega} P_0(x) = 1.
	\end{array}
	\right.
\end{align}

The null distribution $P_0=P_{0b}$ in \re{P0MHtype} involves one single nuisance parameter $\xi=b$. In the special case when beneficial and deleterious mutations have the same rate ($b=1$), this procedure generates a uniform distribution $P_0(x)\equiv 2^{-d}$. On the other hand, states $x$ with many functioning parts will be harder to reach by the Markov process $\tilde{\mathbf \Pi}_0$ when beneficial mutations occur less frequently than deleterious ones ($0<b<1$), resulting in smaller values of $P_0(x)$. The distribution under the alternative hypothesis, $P=P_{\tilde{\theta}bt}$, involves the nuisance parameter $b$, the time point $t$ at which the state of the population is recorded, and $\tilde{\theta}=(a,\theta)$, the two parameters that determine how much background information the MH-type evolutionary algorithm makes use of. For simplicity $a$ and $b$ are here regarded as constants and we only include $\theta$ and $t$ in the notation. This gives rise to an active information
\begin{equation}
I^+(\theta,t) = \log \frac{P_{\theta}\left(X_t=(1,\ldots,1)\right)}{P_{0}\left(X_t=(1,\ldots,1)\right)}.
\label{I+Moran}
\end{equation}

The MH-type algorithm is studied for $d=5$, and illustrated in Figures \ref{Fig:1}-\ref{Fig:3}. Note that the functional information $I_{f0}$ is a decreasing function of $b$, since it is more surprising to find a working molecular machine by chance when the rate of beneficial mutations $b$ is small. Moreover, the active information $I^+(\theta)=\lim_{t\to\infty}I^+(\theta,t)$ for the equilibrium distribution of the Markov chain as well as the active informations $I^+(\theta,t)$ and $I_s^+(\theta,t)$ for a system in non-equilibrium, without and with stopping, are increasing functions of $\theta$, and decreasing functions of $a$ and $b$. The smaller $a$ or $b$ is, the more external information can be infused to increase the probability of reaching the fine-tuned state of a working molecular machine $(1,\ldots,1)$. When $a$ is small, to leave this state once it is reached becomes more difficult, and consequently $I_s^+(\theta,t)$ is only marginally larger than $I(\theta,t)$.    
\end{exa}

\begin{figure}[!t]
	\centering
	\includegraphics[width=2.5in]{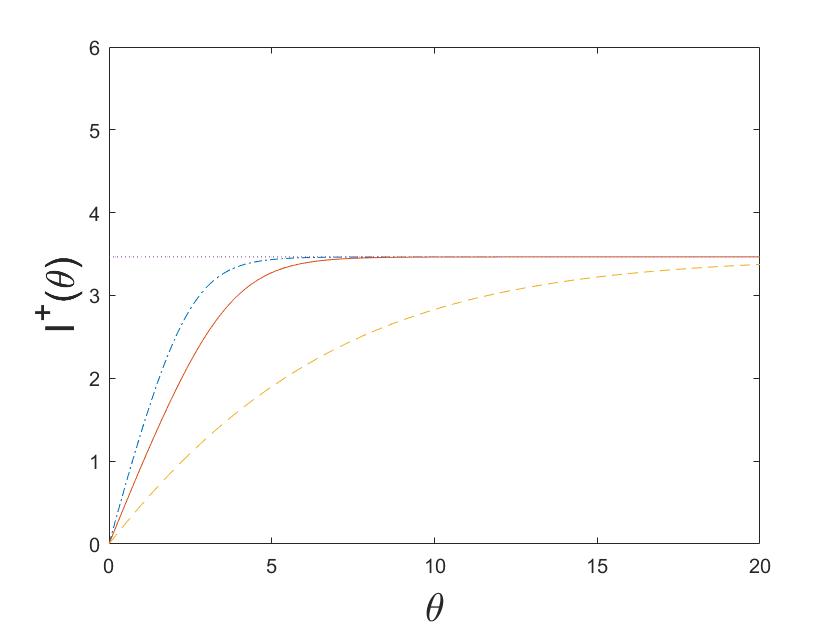}
	\caption{Plot of $I^+(\theta)=\lim_{t\to\infty} I^+(\theta,t)$ in \eqref{I+Moran} as a function of $\theta$ for a system of molecular machines with transition kernel \re{alphaN}, proposal distribution \re{q}, and null distribution \re{P0MHtype}. The system has $d=5$ components, $b=1.0$, and $a=-0.2$ (dash-dotted), $a=0$ (solid) and $a=0.2$ (dashed). The horizontal dotted line corresponds to the functional information $I_{f0}=3.47$.}
	\label{Fig:1}
\end{figure}

\begin{figure}[!t]
	\centering
	\includegraphics[width=2.5in]{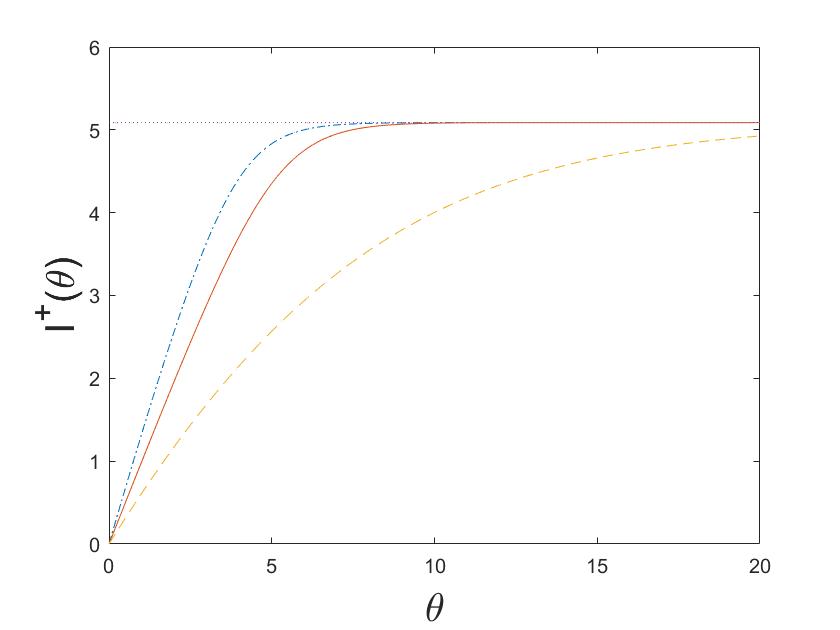}
	\caption{Plot of $I^+(\theta)=\lim_{t\to\infty} I^+(\theta,t)$ in \eqref{I+Moran} as a function of $\theta$ for a system of molecular machines with transition kernel \re{alphaN}, proposal distribution \re{q}, and null distribution \re{P0MHtype}. The system has $d=5$ components, $b=0.5$, and $a=-0.2$ (dash-dotted), $a=0$ (solid) and $a=0.2$ (dashed). The horizontal dotted line corresponds to the functional information $I_{f0}=5.09$.}
	\label{Fig:2}
\end{figure}

\begin{figure}[!t]
	\centering
	{\includegraphics[width=1.5in]{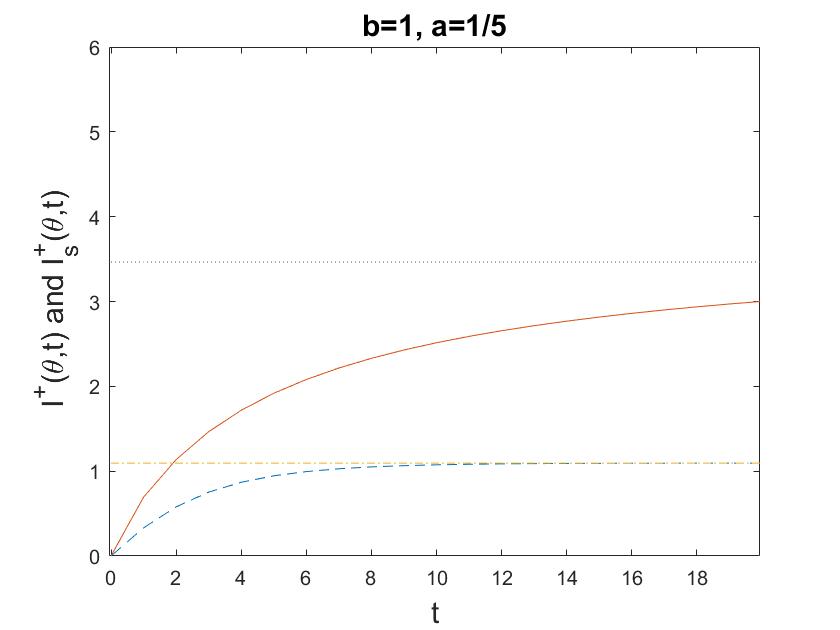}}
	\qquad
	{\includegraphics[width=1.5in]{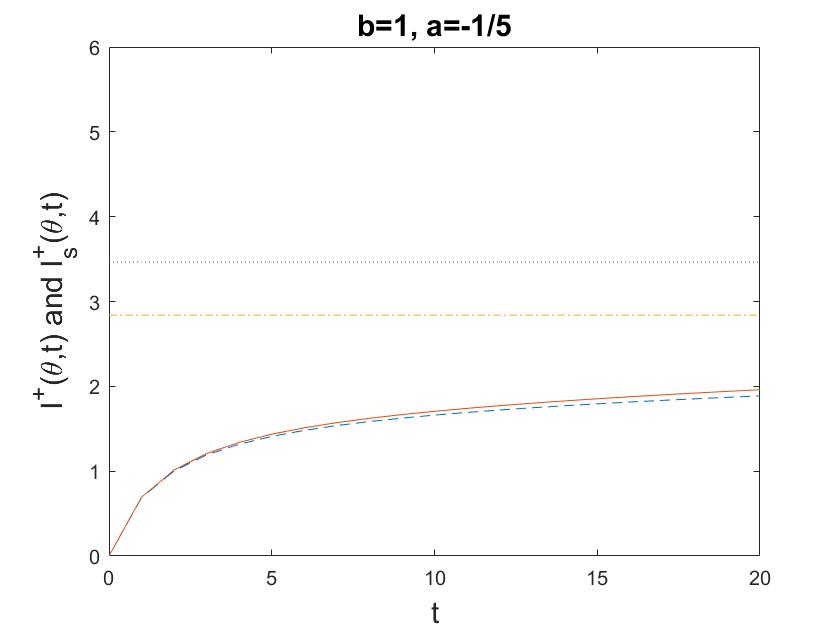}}
	\\
	{\includegraphics[width=1.5in]{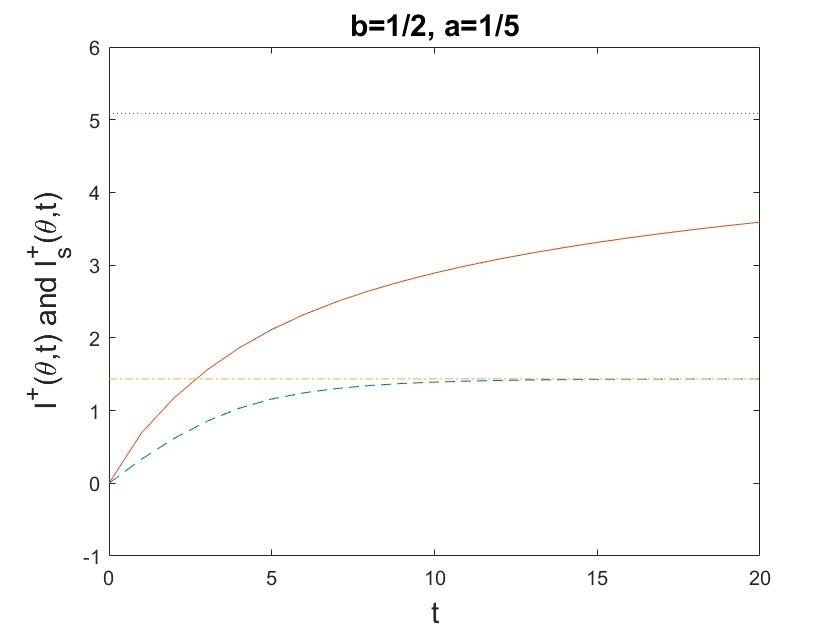}}
	\qquad
	{\includegraphics[width=1.5in]{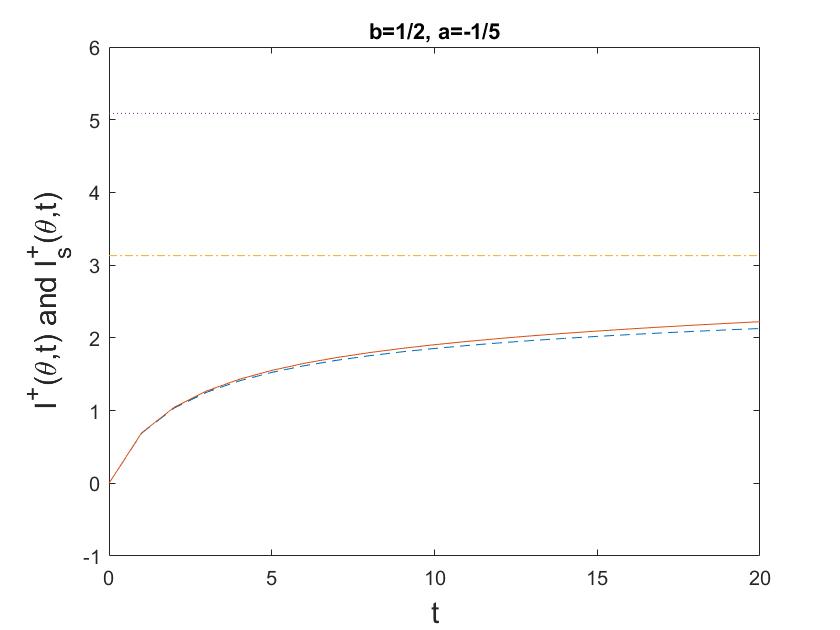}}
	\caption{Plot of $I^+(\theta,t)$ in \eqref{I+Moran} (dashed) and $I^+_s(\theta,t)$ (solid) as a function of $t$ for a system of molecular machines with transition kernel \re{alphaN}, proposal distribution \re{q}, and null distribution \re{P0MHtype}. The system has $d=5$ components and $\theta=2.5$. The upper (lower) row corresponds to $b=1$ ($b=0.5$), whereas the left (right) column corresponds to $a=0.2$ ($a=-0.2$). The horizontal lines 	in each figure illustrate $I^+(\theta)$ (dash-dotted) and the functional information $I_{f0}$ (dotted).}
	\label{Fig:3}
\end{figure}

\begin{exa}[Evolutionary programming algorithms]
Suppose $\Omega=\Omega_{\text{ind}}^N$ is a set of genetic variants from some genomic region, $x=(x_1,\ldots,x_N)$, for the members of a population of size $N$. That is, $x_k\in\Omega_{\text{ind}}$ is the variant of this genomic region for individual $k$. If, for instance, the region codes for the molecular machine of Example \ref{Ex:Moran}, we let $x_k=(x_{k1},\ldots,x_{kd})\in \{0,1\}^d=\Omega_{\text{ind}}$, with $x_{kj}=1$ or 0 depending on whether component $j$ of this machine works or not for individual $k$. Let $g(x_k)$ be the biological fitness, or expected number of offspring, of $k$. In the context of molecular machines, the logarithm of $g(x_k)$ could be a function of the number of functioning parts of a machine of type $x_k$. The specificity function of a population in state $x$ is the average fitness
$$
f(x) = \frac{1}{N}\sum_{k=1}^N g(x_k)
$$
of its individuals. The targeted set $A$ in \re{A} corresponds to all genetic profiles with an average fitness at least $f_0$. This type of model is frequently used in genetic programming as well as in other types of evolutionary programming algorithms to mimic the evolution of $N$ individuals over time \cite{Mitchell1996, Vikhar2016}. Typically, the output $X=X_t$ of the evolutionary algorithm is the last step of a simulation $\{X_s=(X_{s1},\ldots,X_{sN})\}_{s=0}^t$ of the population over $t$ generations. Once the distributions $P_0=P_{0\xi t}$ and $P=P_{\theta\xi t}$ of $X$ are found under the null hypothesis $H_0$ and the alternative hypothesis $H_1$, the actinfo $I^+$ can be computed, according to \eqref{I+}. This actinfo quantifies, on a logarithmic scale, how much more likely it is for the average fitness of the population to exceed $f_0$ at time $t$, for a population with externally infused information ($H_1$) compared to an evolutionary process where no such external information is used ($H_0$). For instance, if a molecular machine needs all its parts in order to function ($g(x_k)=1(|x_k|=d)$), then the actinfo at time $t$ equals
\begin{equation}
I^+ = I^+(\theta,\xi,t) = \log \frac{P_{\theta\xi} \left(|\{k; \, 1\le k \le N, X_{tk}=(1,\ldots,1)\}|\ge f_0 N\right)}
{P_{0\xi} \left(|\{k; \, 1\le k \le N, X_{tk}=(1,\ldots,1)\}|\ge f_0 N\right)}.
\label{I+EvAlg}
\end{equation}

Since the state space $\Omega$ is very large, it is often complicated to find explicit, analytical expressions for the actinfo $I^+$ in \eqref{I+EvAlg}. Suppose the nuisance parameters $\xi$ of the null distribution $P_0=P_{0\xi}$ are known. This makes  the framework of Section \ref{Sec:H0known} applicable, running the evolutionary algorithm $n$ times. That is, $n$ i.i.d.\ copies $\{X_{is}\}_{s=0}^t$ of the population trajectory are generated up to time $t$ for $i=1,\ldots,n$. Then $X_i=X_{it}=(X_{it1},\ldots,X_{itN})$, $i=1,\ldots,n$, are used for computing an estimate $\hI_n^+$ of the actinfo, and test for fine-tuning, according to Section \ref{Sec:H0known}. 

Recall the fixed state assumption of Example \ref{Ex:Moran}, whereby all individuals of the population, at any time point, have the same state. Such an assumption is only realistic when $N\mu \ll 1$, that is, when either the mutation rate $\mu$ and/or the population size $N$ is small. It corresponds to a scenario where $P_0$ and $P$ put all their probability masses along the diagonal
\begin{align}\label{Omegadiag}
\Omega_{\text{diag}} = \{x\in\Omega; \, x_1=\ldots=x_N\}
\end{align}
of $\Omega$. Since \re{Omegadiag} is equivalent to the reduced state space $\Omega_{\text{ind}}$, the fixed state assumption greatly simplifies the analysis. For instance, it often makes it possible to find analytical expressions for the actinfo $I^+$, rather than having to estimate it.    
\end{exa}

\section{Discussion}\label{Sec:Disc}

In this article a general statistical framework is provided for using active information to quantify the amount of pre-specified external knowledge an algorithm makes use of, or equivalently, how tuned the algorithm is. The theory is based on quantifying, for each state $x$, how specified it is by means of a real-valued function $f(x)$. An algorithm with external information either makes use of knowledge of $f$ directly, or at least it incorporates knowledge that tends to move the output of the algorithm towards more specified regions. The Metropolis-Hastings Markov chain incorporates knowledge of $f$ directly, in terms of the acceptance probability of proposed moves. The learning ability of this algorithm was analyzed by studying its active information, with or without stopping, when the targeted set of highly specified states is reached. When independent outcomes of an algorithm are available, nonparametric and parametric estimators of the actinfo of the algorithm were also developed, as well as nonparametric and parametric tests of FT. 

This work can be extended in different ways. A first extension is to find conditions under which the actinfo $I^+(\theta,t)$ of a stochastic algorithm based on a random start (according to the null distribution of a non-guided algorithm) followed by $t$ iterations of the Metropolis-Hastings Markov chain (without stopping) is a non-decreasing function of $t$. We conjecture that this is typically the case but have not obtained any general conditions on the distribution $q$ of proposed candidates for this result to hold. 

A second extension is to widen the notion of specificity, so that not only the functionality $f(x)$ but also the rarity $P_0(x)$ of the outcome $x$ under the null distribution is taken into account. A class of such specificity functions is 
\begin{align}\label{gtheta2}
	g_\theta(x) = \theta f(x) - \log P_0(x),
\end{align}
where $\theta>0$ is a parameter that controls the tradeoff between scenarios where either functionality or rarity under the null is the most important determinant of specificity. The case $\theta=0$ in \eqref{gtheta2} corresponds to function having no impact, so that $g_0(x)$ reduces to Shannon's self information of $x$. The case $g_1(x)$ was proposed in \cite{Montanez2018}, whereas $g_\theta(x)$ is solely determined by $f(x)$ in the limit when $\theta$ gets large.      

A third extension is to generalize the notion of actinfo to include not only the probability of reaching a targeted set of highly specified states $A$ under $H_0$ and $H_1$, but also account for the conditional distribution of the states within $A$, given that $A$ has been reached. This is related to the way in which {\it functional sequence complexity} generalizes functional information \cite{AbelTrevors2005, DurstonChiu2005, DurstonChiu2011, DurstonEtAl2007}. Let $H(Q)=-\sum_x Q(x)\log [Q(x)]$ refer to the Shannon entropy of a distribution $Q$, whereas $H(Q_A)$ is the Shannon entropy of the corresponding conditional distribution $Q_{A}(x) = Q(x|A)$, given that $A$ has been reached. The functional sequence complexity 
\begin{align*}
		\text{FSC}_0 &= H(P_0) - H(P_{0A})\\
		&=  E_{P_0}\left\{\log [P_0(X \mid A)] \mid X\in A\right\} - E_{P_0}\{\log [P_0(X)]\}
\end{align*}
is the reduction in entropy, under the null hypothesis $H_0$ of the highly specified states in $A$, compared to the entropy under $H_0$ of all states in $\Omega$. $\text{FSC}_0$ then reduces to the functional  information $I_{f0}$ when $P_0$ is uniform over $\Omega$. In a similar vein, the {\it active uncertainty reduction} is introduced:
\begin{align*}
	\text{UR}^+ &= \sum_{x\in A} P_A(x) \log P(x) - \sum_{x\in A} P_{0A}(x)\log P_0(x)\\
	&= E_P[\log P(X)|X\in A] - E_{P_0}[\log P_0(X)|X\in A].
\end{align*} 

Then $\text{UR}^+=I^+$ when $P_{0A}$ and $P_A$ are uniformly distributed on $A$. This happens, for instance, when $P_0$ has a uniform distribution on $\Omega$ and $P=P_\theta$ for some $\theta>0$, and if \eqref{f0fmax} holds. The properties of $\mbox{UR}^+$ deserve to be analyzed in more detail, for instance investigate how it differs from the actinfo $I^+$. 

A fourth extension would be to apply the concept of actinfo to other genetic models. For instance, Example \ref{Ex:Moran} is the first time that, to our knowledge, actinfo is applied to the Moran model. In the past though, actinfo was used in population genetics to study fixation times for the Wright-Fisher model of population genetics, a model for which time is discrete and generations do not overlap \cite{DiazMarks2020b}.

\section{Proofs}\label{Sec:Proofs}

{\bf Proof of Proposition \ref{Prop:theta}.} Introduce
	\begin{equation}
		\begin{aligned}\label{JK}
			J(\theta) &= \sum_{x\in A^c} \exp\{\theta[f(x)-f(x_0)]\}P_0(x),\\
			K(\theta) &= \sum_{x\in A} \exp\{\theta[f(x)-f(x_0)]\}P_0(x),
		\end{aligned}
	\end{equation}
	when $\Omega$ is finite, and replace the sums in \eqref{JK} by integrals when $\Omega$ is continuous. Then 
	\begin{align}\label{PJK}
		P_\theta(A) &=  \exp[\theta f(x_0)]K(\theta) / \{\exp(\theta f(x_0))[J(\theta)+K(\theta)]\} \nonumber \\
		&= K(\theta) / [J(\theta) + K(\theta)] \\
		&= 1 / [J(\theta)/K(\theta) +1]. \nonumber
	\end{align}

	Since $P_0(A)<1$, it follows that $J(\theta)$ is a strictly decreasing function of $\theta \ge 0$, whereas $K(\theta)$ is a non-decreasing function of $\theta$. 		From this, it follows that $P_\theta(A)$ is a strictly increasing function of $\theta$, and consequently $I^+(\theta) = \log[P_\theta(A)/P_0(A)]$ is a strictly 			increasing function of $\theta$ as well. 

	Moreover, $K(\theta)\ge P_0(A)>0$ for all $\theta\ge 0$, and $J(\theta)\to 0$ as $\theta\to\infty$ follows by dominated convergence. In conjunction with 			\eqref{PJK} this implies $P_\theta(A)\to 1$ and $I^+(\theta)\to I_{f0}$ as $\theta\to\infty$. 

\par\bigskip

{\parindent=0pt {\bf Proof of Proposition \ref{Prop:thetat}} Equation \eqref{I+thetat} follows from \eqref{bP0}, \eqref{bPiConv} and the fact that } 
	\begin{align*}
		P_0(A) &= \sum_{x\in A} P_0(x) = \mathbf P_0\mathbf v,\\
		P_{\theta t}(A) &= \sum_{x\in A} P_{\theta t}(x) = \mathbf P_{\theta t}\mathbf v = \mathbf P_0 \mathbf \Pi_\theta^t \mathbf v,
	\end{align*}
	since $\mathbf v$ is a column vector of length $|\Omega|$ with ones in positions $x\in A$ and zeros in positions $x\in A^c$. 

	Equation \eqref{I+Conv} is equivalent to proving that 
	\begin{align*}
		P_{\theta t}(A) \to P_\theta(A) \text{ as }t\to\infty.
	\end{align*}
	But this follows from the fact that $P_\theta$ is the equilibrium distribution of the Markov chain with transition kernel \eqref{bPi}. That is, letting $t\to\infty$ in 		\eqref{bPiConv} we find that
	\begin{align*}
		\mathbf P_{\theta t} = \mathbf P_0\mathbf \Pi_\theta^t\to\mathbf P_\theta, 
	\end{align*} 
	and therefore
	$$
		P_{\theta t}(A) = \mathbf P_{\theta t}\mathbf v \to \mathbf P_\theta\mathbf v = P_\theta(A), \text{ as } t \to \infty.
	$$

\par\bigskip

{\parindent=0pt {\bf Proof of Proposition \ref{Prop:thetats}} Equation \eqref{I+Ineq} follows from the definitions of $I^+(\theta,t)$ and $I_s^+(\theta,t)$ in \eqref{I+thetat} and \eqref{I+thetats}, and the fact that} 
	$$
		P_{\theta t}(A) = P(X_t\in A) \le P(X_{t\wedge T}\in A) = P_{\theta ts}(A),
	$$
	where the inequality is a consequence of the definition of $T$ in \eqref{T}. Since 
	$$
		P_{\theta ts}(A) = P(T\le t) \le P(T\le t+1) = P_{\theta , t+1,s}(A),
	$$
	we have proved that $I_s^+(\theta,t)$ is non-decreasing in $t$. Equation \eqref{I+Lim} follows from the definition of $I_s^+(\theta,t)$ and the fact that 
	\begin{equation}\label{TLim}
		\lim_{t\to\infty} P_{\theta ts}(A) = P(T<\infty) = 1.
	\end{equation}
	The last equality of \eqref{TLim} is a consequence of the fact that the Markov chain with transition kernel $\mathbf \Pi_\theta$ is irreducible, so that any state 		$x\in\Omega$ will be reached with probability 1. In particular, the targeted set $A$ will be reached with probability 1. In order to verify \eqref{I+ET}, we first 		deduce 
	$$
		P(T>t) = 1 - P_0(A)e^{I^+_s(\theta,t)}
	$$
	from \eqref{Ptheta ts2}, and then we make use of the equality 
	$$
		E(T) = \sum_{t=0}^\infty P(T>t).
	$$

\par\bigskip

{\parindent=0pt {\bf Proof of Proposition \ref{Prop:LD}.} Since $n\hat Q(A)\sim \text{Bin}(n,Q(A))$ has a binomial distribution, it follows from the Central Limit Theorem that}
	\begin{equation}\label{CLT}
		\sqrt{n}(\hat Q(A)-Q(A)) \stackrel{\cal L} \longrightarrow N(0,Q(A)[1-Q(A)]),
	\end{equation}
	as $n\to\infty$. Notice that $\hat I^+ = g(\hat Q(A))$, where $g(Q)=\log [Q/P_0(A)]$ and $g^\prime(Q)=1/Q$. Equation \eqref{AsNNP} follows from the Delta 		Method (see, e.g., Theorem 8.12 of \cite{LehmannCasella1998}) and the fact that
	$$
		V = g^\prime(Q(A))^2 \cdot Q(A)[1-Q(A)].
	$$

	In order to establish \eqref{LD}, to begin with, it follows from \eqref{hI+} and the definition of $p_{\text{min}}$ that
	\begin{align*}
		P_{H_0}(\hat I^+ \ge I_\text{min}) &= P_{H_0}(\hat Q(A) \ge p_\text{min})\\ 
		&= P_{H_0}\left(\frac{1}{n}\sum_{i=1}^n Y_i \ge p_\text{min}\right),
	\end{align*}
	where $Y_i=I(X_i\in A) \sim \mbox{Be}(p_0)$ are independent Bernoulli variables under $H_0$ with success probability $p_0=P_0(A)$. It follows from Large 		Deviations theory that \eqref{LD} holds, with 
	\begin{equation}\label{CLF}
		C = \sup_{\phi >0} [\phi p_\text{min} - \lambda(\phi)]
	\end{equation}
	the Legendre-Fenchel transformation, and 
	\begin{equation}\label{la}
		\lambda(\phi) = \log E[\exp(\phi Y)] = \log [1+p_0 (e^\phi-1)]
	\end{equation}
	the cumulant generating function of $Y$ \cite[pp.~529-533]{Kallenberg2021b}. Inserting \eqref{la} into \eqref{CLF} it can be seen that the maximum in \eqref{CLF} is given by \eqref{C}. 

\par\bigskip

{\parindent=0pt {\bf Proof of Proposition \ref{Prop:Param}.} 	In order to verify \eqref{AsNP}, we will first show that the estimator \eqref{hth} of the tilting parameter $\theta$ is asymptotically normal} 
	\begin{equation}\label{AsN}
		\sqrt n(\hat \theta_n-\theta^\ast) \stackrel{\cal L} \longrightarrow N(0,U)\mbox{ as }n\to\infty,
	\end{equation}
	with asymptotic variance
	\begin{equation}\label{V}
		U = \frac{\text{Var}_Q[f(X)]}{\text{Var}^{\,\, 2}_{P_{\theta^\ast}}[f(X)]}.
	\end{equation}
	
	To this end, let ${}^\prime$ refer to derivatives with respect to the tilting parameter $\theta$. Define the score function 
	$$
		\psi_\theta(x) = \frac{d\log P_\theta(x)}{d\theta} = \frac{P_\theta^\prime(x)}{P_\theta(x)}
	$$
	and its derivative 
	$$
		\psi_\theta^\prime(x) = \frac{d\psi_\theta(x)}{d\theta}.
	$$
	
	It is a standard result from the asymptotic theory of maximum likelihood estimation and $M$-estimation (see, e.g.,~Chapter 6 of \cite{LehmannCasella1998}) that \eqref{AsN} holds with asymptotic variance
	\begin{equation}\label{V2a}
		U = \frac{\text{Var}_Q[\psi_{\theta^\ast}(X)]}{E^2_Q[\psi^\prime_{\theta^\ast}(X)]}.
	\end{equation}
	
	To simplify \eqref{V2a}, notice that the score function can be written as 
	\begin{equation}\label{psith}
		\psi_\theta(x) = f(x) - \frac{M^\prime(\theta)}{M(\theta)} = f(x) - E_{P_\theta}[f(X)]
	\end{equation}
	for the exponential family of tilted distributions \eqref{Ptheta}-\eqref{Mtheta}. From this it follows that
	$$
		\psi^\prime_\theta(x) = \frac{M^{\prime\prime}(\theta)}{M(\theta)} - \left( \frac{M^\prime(\theta)}{M(\theta)}\right)^2
		= \text{Var}_{P_\theta}[f(X)]
	$$
	is a constant, not depending on $x$. Inserting the last two displayed equations into \eqref{V2a}, the formula in \eqref{V} for the asymptotic variance of $\hat 		\theta$ is obtained.  As a next step we notice that 
	\begin{equation}\label{hI+g}
		\hat I^+ = g(\hat \theta),
	\end{equation}
	where 
	\begin{equation}\label{gtheta}
		g(\theta) = \log \frac{P_\theta(A)}{P_0(A)} = \log h(\theta) - \log P_0(A),
	\end{equation}
	and 
	\begin{equation}\label{h}
		h(\theta) = P_\theta(A) = \frac{\sum_{x\in A} e^{\theta f(x)} P_0(x)dx}{M(\theta)}
	\end{equation}
	follows from the definition of $P_\theta(x)$ in \eqref{Ptheta}. 
	
	Differentiating \eqref{h} with respect to $\theta$, we find that
	\begin{align}\label{hprime}
	      \begin{aligned}
		h^\prime(\theta) &= \sum_{x\in A} f(x)e^{\theta f(x)} P_0(x)dx /M(\theta)\\ 
		&- M^\prime(\theta)\sum_{x\in A} e^{\theta f(x)} P_0(x)dx /M^2(\theta).
		\end{aligned}
	\end{align}
	And it follows from the RHS of \eqref{hprime} that
	\begin{align}
		\begin{aligned}
		h^\prime(\theta) &= E_{P_\theta}[f(X)I(f(X)\ge f_0)] - P_\theta(A) E_{P_\theta}[f(X)]\\
					&= \text{Cov}_{P_\theta}[f(X),I(f(X)\ge f_0)].
		\end{aligned}
	\end{align}

	Then we combine \eqref{gtheta} and \eqref{hprime}, and obtain
	\begin{equation}\label{gprimea}
		g^\prime(\theta) = \frac{h^\prime(\theta)}{h(\theta)} = \frac{ \text{Cov}_{P_\theta}[f(X),I(f(X)\ge f_0)]}{P_\theta(A)}.
	\end{equation}
	Finally we use the Delta Method to conclude that $\hat I^+$ is an asymptotic normal estimator \eqref{AsNNP} of $I^+(\theta^\ast)$, with asymptotic variance $V 	= g^\prime(\theta^\ast)^2 U$, which, in view of \eqref{V} and \eqref{gprimea}, agrees with \eqref{VPar}.  

	In order to prove the large deviation result \eqref{LD2} for the parametric test of FT, let $\theta_\text{min}$ be the value of the tilting parameter that 	satisfies $P_{\theta_\text{min}}(A)=p_\text{min}=P_0(A)\exp(I_{\text{min}})$. Then notice that 
	\begin{align*}
			P_{H_0}(\hat I^+ \ge I_\text{min}) &= P_{H_0}(\hat Q(A) \ge p_\text{min})\\
								&= P_{H_0}(\hat \theta \ge \theta_\text{min})\\
								&= P_{H_0}(\sum_{i=1}^n \psi_{\theta_\text{min}}(X_i)/n \ge 0)\\
								&= P_{H_0}\left(\sum_{i=1}^n f(X_i)/n \ge E_{p_\text{min}}[f(X)]\right), 
	\end{align*}
	 where in the third step we utilized that $\hat \theta \ge \theta_\text{min}$ is equivalent to the derivative of the log likelihood of data being non-negative at 			$\theta_\text{min}$, and in the fourth step we made use of \eqref{psith} and introduced $p_\text{min} = P_{\theta_\text{min}}$. But this last line is a large 			deviations probability. It then follows from large deviations theory that \eqref{LD2} holds, with $C$ the Legendre-Fenchel transformation in \eqref{C2}.   

\par\bigskip

{\parindent=0pt {\bf Proof of Proposition \ref{Prop:Nuisance}.} Since the bias corrected empirical actinfo 
\begin{equation}
\hI_n^+ - B = \log \frac{\hQ(A)}{P_{0\xi}(A)}
\lb{hI+B}
\end{equation}
behaves like \re{hI+}, with $P_0=P_{0\xi}$, the asymptotic normality result for the nonparametric version of the estimator of $I_Q^+$ follows from Proposition \ref{Prop:LD}. 
\par\medskip
For the parametric version of the estimator of $I_Q^+$ we will (briefly) generalize the asymptotic normality proof of Proposition \ref{Prop:Param}. It follows from \re{hI+max} and \re{hQA3} that
$$
\hI_n^+ = g(\hth,\hxi),
$$
where
\begin{equation}
g(\theta,\xi) = \log \frac{P_{\theta\xi}(A)}{P_{0\text{max}}(A)}.
\lb{g}
\end{equation}
Making use of the delta method, it follows that the asymptotic variance of the parametric version of $\hI_n^+$ equals 
\begin{equation}
V = g^\prime(\theta^\ast,\xi^\ast) \mbox{AsVar}(\hth,\hxi)  g^\prime(\theta^\ast,\xi^\ast)^T,
\lb{VProof}
\end{equation}
with the asymptotic variance of $(\hth,\hxi)$ defined through
$$
\sqrt{n}\left( (\hth,\hxi)-(\theta^\ast,\xi^\ast)\right) \Lto N(0,\mbox{AsVar}(\hth,\hxi))
$$
as $n\to\infty$. Since $(\hth,\hxi)$ in \re{hthhxi} is an $M$-estimator, it follows that its asymptotic variance equals 
\begin{equation}
\mbox{AsVar}(\hth,\hxi) = E[\psi_{\theta^\ast \xi^\ast}^\prime (X)]^{-1} E[\psi^T_{\theta^\ast \xi^\ast}(X)\psi_{\theta^\ast \xi^\ast}(X)] E[(\psi_{\theta^\ast \xi^\ast}^\prime)^T(X)]^{-1}.
\lb{AsVar}
\end{equation}
The gradient of \re{g} is
\begin{equation}
g^\prime(\theta,\xi) = \frac{P_{\theta\xi}^\prime(A)}{P_{\theta\xi}(A)} = E[\psi_{\theta\xi}(X)|X\in A], 
\lb{gprime}
\end{equation}
where $\psi_{\theta\xi}= P_{\theta\xi}^\prime(x) / P_{\theta\xi}(x)$ is the likelihood score function for the combined para\-meter vector $(\theta,\xi)$. Putting things together, the asympotic variance formula \re{VPar2} for the parametric version of $\hI_n^+$ follows from \re{VProof}-\re{gprime}.    
\par\medskip
The significance level of the FT test can be written as 
$$
P_{0\xi}(\hI_n^+ \ge I_{\text{min}}) = P_{0\xi}(\hI_n^+ - B \ge I_{\text{min}}-B).
$$
Since $p_{\text{min}} = P_{0\xi}(A)\exp(I_{\text{min}})$, we have that
\begin{equation}
I_{\text{min}}-B =   \log \frac{p_{\text{min}}e^{-B}}{P_{0\xi}(A)}.
\lb{IminB}
\end{equation}
From this and \re{hI+B} it follows that the nonparametric test of FT behaves as the correspon\-ding nonparametric test of Proposition \ref{Prop:LD}, with the null probability $P_0(A)$ replaced by $P_{0\xi}(A)$, and $p_{\text{min}}$ replaced by  $p_{\text{min}}e^{-B}$. Therefore, the large deviation result \re{CNonpar} follows from \re{C}. In a similar way, the large deviation result for the parametric version of the FT-test (in the special case when $\theta$ is a scalar exponential tilting parameter) follows from \re{hI+B}, \re{IminB} and Proposition \ref{Prop:Param}.  

\par\bigskip

{\parindent=0pt {\bf Proof of Proposition \ref{Prop:TwoSamples}.} Because of \re{I+nuisance} and \re{hITwoSamp} we have that
\begin{equation}
\sqrt{n}(\hI_{nn_0}^+-I_Q^+) = \sqrt{n}\log \frac{\hQ(A)}{Q(A)} - \sqrt{\frac{n}{n_0}} \sqrt{n_0}\log \frac{P_{0\hxi}(A)}{P_{0\xi}(A)},
\lb{hInn0Split}
\end{equation}
where
\begin{equation}
\sqrt{n} \log \frac{\hQ(A)}{Q(A)} \Lto N(0,V_1) \mbox{ as }n\to\infty
\lb{V1}
\end{equation}
and
\begin{equation}
\sqrt{n_0} \log \frac{P_{0\hxi}(A)}{P_{0\xi}(A)} \Lto N(0,V_2)\mbox{ as }n_0\to\infty
\lb{V2}
\end{equation}
respectively. It follows from the proofs of Propositions \ref{Prop:LD}-\ref{Prop:Param} that the asymptotic variance for $V_1$ in \re{V1} is the same as $V$ in \re{VNP} and \re{VPar2}, for the nonparametric and prametric versions of $\hQ(A)$ respectively. 
The asymptotic variance $V_2$ in \re{V2} is given by \re{V2TwoSamp}. This is proved using the delta method (similarly as for Proposition \ref{Prop:Nuisance}), making use of the fact that $\hxi$ is the maximum likelihood estimator of $\xi$ with asymptotic variance that is the inverse $E[\psi_\xi^T(X)\psi_\xi(X)]^{-1}$ of the Fisher information matrix. The asymptotic normality result \re{hIAsNTwoSamp} then follows from \re{hInn0Split}-\re{V2}, the fact that $n/n_0\to \la$, and the independence of the two samples. 
\par\medskip
The large deviations results are proved in a similar way as in Proposition \ref{Prop:Nuisance}, replacing $P_{0\text{max}}(A)$ by $P_{0\hxi}(A)$. Using statistical consistency $\hxi\pto \xi$ as $n_0\to\infty$, it follows that the large deviation rates $C$ of Proposition \ref{Prop:TwoSamples}, for the nonparametric and parametric versions of the FT tests, are the same as in Proposition \ref{Prop:Nuisance}, with bias term $B=0$. 

\par\bigskip

{\parindent=0pt {\bf Details from Example \ref{Ex:Moran}.} In order to prove that the Metropolis-Hastings type Markov chain \eqref{pi} with acceptance probabilities \eqref{alphaN} has equilibrium distribution $P_\theta$, we first notice that for any pair of states $x\ne y$, the flow of probability mass
\begin{align}\label{Flow}
		&P_\theta(x)\pi_{\theta}(x,y) \nonumber \\
		&=P_\theta(x)q(x,y)\alpha_{\theta}(x,y) \nonumber \\
		&=\frac{P_0(x)e^{\theta f(x)}}{M(\theta)} q(x,y)\cdot  C \left[ \frac{e^{\theta f(y)} P_0(y)q(y,x)}{e^{\theta f(x)} P_0(x)q(x,y)} \right]^{1/2} \nonumber \\
		&=C \frac{\left(  e^{\theta f(x)} P_0(x)q(x,y) e^{\theta f(y)} P_0(y) q(y,x)\right)^{1/2}}{M(\theta)}
\end{align}
from $x$ to $y$ is symmetric with respect to $x$ and $y$. Therefore, the flow $P_\theta(y)\pi_{\theta}(y,x)$ of probability mass in the opposite direction, from $y$ to $x$, is the same as in \eqref{Flow}. A Markov chain with this property is called {\it reversible} \cite[pp.~11-12]{Popov2021}. But it is well known that $P_\theta$ is a stationary distribution if the Markov chain is reversible with reversible measure $P_\theta$ \cite[p.~238]{GrimmettStirzaker2001}. If, additionally, the proposal distribution $q$ is such that it is possible to move between any pair of states in a finite number of steps, it follows that the Markov chain is irreducible and hence that $P_\theta$ is its unique stationary distribution, which is also the equilibrium distribution of the Markov chain \cite[p.~232]{GrimmettStirzaker2001}. 

Next we will motivate formula \eqref{alphaMoran} for the acceptance probability of a Moran model. Assume that the population evolves over time as a Moran model, and that all individuals have type $x$. If one individual mutates from $x$ to $y$, because of \eqref{sx}, the relative fitness between the $N-1$ individuals of type $x$ and the newly mutated individual of type $y$ is 
\begin{equation}\label{s}
	s = \frac{e^{\theta f(y)/N}}{e^{\theta f(x)/N}} = e^{\theta[f(y)-f(x)]/N}. 
\end{equation}

From the theory of Moran models (e.g., \cite{HossjerBechlyGauger2021,KomarovaEtAl2003}), it is well known that the fixation probability for the newly mutated individual is 
\begin{equation}\label{beta}
	\beta_N(s) = \left\{
	\begin{array}{ll}
		(1-s^{-1})/(1-s^{-N}), & s \ne 1,\\
		1/N , & s=1.
	\end{array}\right.
\end{equation}

Inserting \eqref{s} into \eqref{beta} we find (when $s\ne 1$, or equivalently when $\Delta=\theta[f(y)-f(x)]\ne 0$), that
\begin{equation*}
	\beta_N(s) = \frac{1-e^{-\Delta/N}}{1-e^{-\Delta}} \approx \frac{1}{N} \cdot \frac{\Delta}{1-e^{-\Delta}} \approx
	\frac{1}{N} \cdot (1 + \frac{\Delta}{2}),
\end{equation*}
which is equivalent to \eqref{alphaMoran}. 

\begin{adjustwidth}{-\extralength}{0cm}

\reftitle{References}

\bibliography{daangapaBibliography.bib}

\end{adjustwidth}

\end{document}